\documentclass[mathpazo]{cicp}

\begin{document}
%%%%% title : short title may not be used but TITLE is required.
% \title{TITLE}
% \title[short title]{TITLE}
\title{A GPU-accelerated Cartesian grid method for PDEs on irregular domain}

%%%%% author(s) :
% single author:
% \author[name in running head]{AUTHOR\corrauth}
% [name in running head] is NOT OPTIONAL, it is a MUST.
% Use \corrauth to indicate the corresponding author.
% Use \email to provide email address of author.
% \footnote and \thanks are not used in the heading section.
% Another acknowlegments/support of grants, state in Acknowledgments section
% \section*{Acknowledgments}
% \author[O.~Author]{Only Author\corrauth}
% \address{School of Mathematical Sciences, Beijing Normal University,
% Beijing 100875, P.R. China}
% \email{{\tt author@email} (O.~Author)}

% multiple authors:
% Note the use of \affil and \affilnum to link names and addresses.
% The author for correspondence is marked by \corrauth.
% use \emails to provide email addresses of authors
% e.g. below example has 3 authors, first author is also the corresponding
%      author, author 1 and 3 having the same address.
\author[Liwei Tan et.~al.]{Liwei Tan\affil{1}\footnotemark[2],
      Minsheng Huang\affil{1}\footnote[2]{These authors contributed equally to this work.}, and Wenjun Ying\affil{2}\comma\corrauth}
\address{\affilnum{1}\ School of Mathematical Sciences, 
            Shanghai Jiao Tong University, 
            Shanghai 200240, P.R. China. \\
          \affilnum{2}\ School of Mathematical Sciences, MOE-LSC and Institute of Natural Sciences, 
          Shanghai Jiao Tong University, Minhang, 
          Shanghai 200240, P.R. China.}
\emails{{\tt wying@sjtu.edu.cn} (W.~Ying)}

% \author[Zhang Z R et.~al.]{Zhengru Zhang\affil{1}\comma\corrauth,
%       Author Chan\affil{2}, and Author Zhao\affil{1}}
% \address{\affilnum{1}\ School of Mathematical Sciences,
%          Beijing Normal University,
%          Beijing 100875, P.R. China. \\
%           \affilnum{2}\ Department of Mathematics,
%           Hong Kong Baptist University, Hong Kong SAR}
% \emails{{\tt zhang@email} (Z.~Zhang), {\tt chan@email} (A.~Chan),
%          {\tt zhao@email} (A.~Zhao)}
% \footnote and \thanks are not used in the heading section.
% Another acknowlegments/support of grants, state in Acknowledgments section
% \section*{Acknowledgments}

%%%%% Begin Abstract %%%%%%%%%%%
\begin{abstract}
   % The kernel-free boundary integral (KFBI) method has successfully solved partial differential equations (PDEs) on irregular domains. All present KFBI methods are implemented on the CPU platform, and this paper presents the algorithms of single-GPU and multiple-GPUs. The KFBI is a Cartesian grid method, which is inherently suitable for GPU parallel computing. On a single GPU, assigning individual threads can control correction, interpolation, and jump calculations. To facilitate the computation of larger-scale problems,  we extend the algorithm to multiple GPUs. Unlike a single GPU, we use the arrowhead decomposition method to solve the interface problem, achieving optimal computing efficiency and load balancing. Numerical examples show that the proposed algorithm is second-order accurate and efficient. Single-GPU solver speeds 50-200 times than traditional CPU while the eight GPUs distributed solver yields up to $60\%$ parallel efficiency.

  The kernel-free boundary integral (KFBI) method has successfully solved partial differential equations (PDEs) on irregular domains. Diverging from traditional boundary integral methods, the computation of boundary integrals in KFBI is executed through the resolution of equivalent simple interface problems on Cartesian grids, utilizing fast algorithms. While existing implementations of KFBI methods predominantly utilize CPU platforms, GPU architecture's superior computational capabilities and extensive memory bandwidth offer an efficient resolution to computational bottlenecks. This paper delineates the algorithms adapted for both single-GPU and multiple-GPU applications. On a single GPU, assigning individual threads can control correction, interpolation, and jump calculations. The algorithm is expanded to multiple GPUs to enhance the processing of larger-scale problems. The arrowhead decomposition method is employed in multiple-GPU settings, ensuring optimal computational efficiency and load balancing. Numerical examples show that the proposed algorithm is second-order accurate and efficient. Single-GPU solver speeds 50-200 times than traditional CPU while the eight GPUs distributed solver yields up to $60\%$ parallel efficiency.
\end{abstract}
%%%%% end %%%%%%%%%%%

%%%%% AMS/PACs/Keywords %%%%%%%%%%%
%\pac{}

%%%% maketitle %%%%%
\maketitle

%%%% Start %%%%%%
\section{Introduction}

Graphics Processing Units (GPU) are co-processors originally devoted to accelerate graphics processing. In the last years, they are extensively used as massively parallel platforms to run general-purpose programs. This practice is mostly known as General-Purpose computing on Graphics Processing Units (GPGPU). This growing trend is confirmed by the number of computers in the top500 ranking that are provided of GPUs, which on November 2023 was 186\cite{Top500}.

One of the areas taking advantage of the capabilities of this kind of accelerators is scientific computing. There are many recent publications describing works that successfully port code from CPU to GPU, achieving important speedups\cite{Navarro2014,Hwu2012}. Elliptic type problems are widely applied in the fields of electrochemistry\cite{QIAN2021109908,DING2019108864}, electromagnetism\cite{Chai2023}, computational fluid dynamics\cite{Greengard1998318, Quartapelle1993NumericalSO}, shape optimisation problems\cite{ZHU2011752,gong2023} and other areas in science\cite{Chapko199747,ZHOU20061,CHENG2006616,SUN2014445},  Solving these problems often requires a substantial computational cost\cite{GarciaRisueno2014}.

An effective and accurate approach for solving elliptical equations is the Kernel-Free Boundary Integral (KFBI) method \cite{ying2007kernel, ying2013kernel, ying2014kernel}, which originates from boundary integral methods. Unlike traditional boundary integral approaches, the KFBI method embeds complex domains into larger, regular computational areas (such as square regions), which are subsequently partitioned using Cartesian grids. The KFBI method not only benefits from the well-conditioning property of the boundary integral equation(BIE) but also avoids explicitly calculating Green's function directly, which is challenging in complex domains\cite{xie2019fourth,ying2013kernel}.
In recent years, the KFBI method has been extensively applied\cite{xie2019fourth,xie2021cartesian,ZHAO2023116163,dong2023kernelfree,zhou2023adi}.

The substantial memory bandwidth and abundant cores in GPU architecture enable the concurrent execution of thousands of computational tasks, leading to significant acceleration. This renders it an efficient solution for addressing computing bottlenecks. More importantly, the GPU architecture suits the Cartesian grid method since each thread easily controls one grid node. Several related works have addressed the GPU acceleration of Cartesian grid methods in the last ten years\cite{REDDY2015287,cuIBM2018Layton, LIANG2014156Solving, huang2015implementation}: the GPU-accelerated VOF by Rajesh Reddy and R. Banerjee\cite{REDDY2015287}, the CUDA-Based IB method by S. K. Layton, A. Krishnan and L. A. Barba\cite{cuIBM2018Layton}, the TVD Runge–Kutta method on multiple GPUs by Liang. S,  Liu. W and Yuan. L\cite{LIANG2014156Solving}, the multiple-GPU based lattice Boltzmann algorithm by Huang. C, Shi. B, He. N and Chai. Z\cite{huang2015implementation}.

As a Cartesian grid method, the essential procedure of the KFBI method involves the correction of irregular points and control points on the interface individually, making it inherently well-suited for GPU-accelerated parallel processing. Furthermore, the KFBI method utilizes an FFT-based solver, well-documented in literature for its suitability with GPU or GPU clusters\cite{Volkov2011,Govindaraju2008,Nandapalan2012,Chen2010,Nukada2012}, to enhance the efficiency of interface problem computations in iterative procedures. In fact, due to the simple grid topology on Cartesian grids, building a highly parallel GPU-accelerated Cartesian grid solver based on the KFBI method is straightforward. The implementation details of the KFBI solver for a single-GPU version are concisely delineated in section $\ref{oneGPU}$, with the corresponding numerical results presented in section $\ref{result}$.

A significant limitation in single-GPU computation is its available memory, which leads to a bottleneck in the size of the computational mesh. In order to expand the calculation scale and improve efficiency, we also study the multiple-GPU architecture in a single node (in one computer), which contains a two-level parallelization: the coarse-grained level composed by GPUs across multiple CPUs at the cost of coordinating GPU-GPU communication via MPI and the fine-grained level formed by CUDA cores on each GPU. Based on the characteristics mentioned above, we have devised a distributed KFBI algorithm that evenly distributes data to each GPU, maximizing the utilization of multiple-GPU parallel capabilities, ensuring computational load balancing, and minimizing inter-GPU communication overhead.

The remainder of the paper is organized as follows. We first introduce  the boundary integral method and the KFBI methods in section $\ref{KFBI}$. Section $\ref{oneGPU}$ describes implementing the KFBI method on a single GPU. The algorithm is then extended to multiple GPUs and summarised in section $\ref{multi_GPU}$. The numerical results are presented in section $\ref{result}$. The advantages, limitations, and prospects for the GPU-based KFBI method are discussed in the final section. 

\section{The kernel-free boundary intergral method} \label{KFBI}
% \subsection{Introduction of KFBI-CPU method}
Suppose $\Omega$ is a bounded irregular and complex domain in $R^2$ or $R^3$ whose boundary $\Gamma = \partial \Omega$ is at least twice continuously differentiable. Let $u(\mathbf{x})$ be an unknown function of $\mathbf{x} \in \mathbf{R}^{d}(d=2, or ~3)$. Assuming $g_{D}(\mathbf{x})$ and $f(\mathbf{x})$ are known function of $\mathbf{x}$  with sufficient smoothness. $\partial_{\mathbf{n}}u(\mathbf{x})$ denotes the normal derivative of $u(\mathbf{x})$ on the boundary, where $\mathbf{n}$ denotes the unit outward normal on $\Gamma$. For simplicity of description, we introduce the KFBI method for the modified Helmholtz equation subject to the Dirichlet boundary condition. 

Consider the modified Helmholtz equation 
\begin{equation}
    \Delta u(\mathbf{x})-\kappa u(\mathbf{x})=f(\mathbf{x}), \quad \text { in } \Omega,
    \label{one_GPU:modified_helmholtz}
\end{equation}
subject to Dirichlet boundary condition 
\begin{equation}
    u(\mathbf{x}) = g_{D}(\mathbf{x}), \quad \text{ on }\Gamma.
    \label{one_GPU:dirichlet_boundary}
\end{equation}

Here, $\kappa$ is assumed to be a positive constant for the modified Helmholtz equation in this paper by default.

\subsection{Boundary integral equation}
\begin{figure}[htpt!]
    \centering
    \includegraphics[width = 0.8\linewidth]{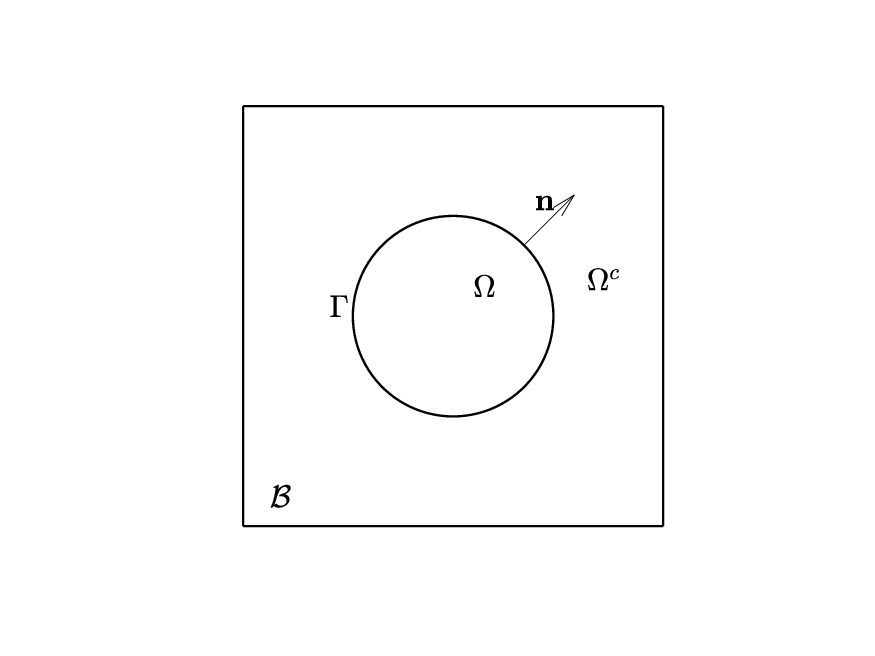}
    \caption{KFBI computation domain}
    \label{kfbi domain}
\end{figure}

As shown in Fig.\,$\ref{kfbi domain}$, to solve the boundary value problem above by the KFBI method, we first embed the irregular domain $\Omega$ into a larger rectangle domain $\mathcal{B} = \Omega \cup \Omega^{c}$. 

According to the standard BIM\cite{aliabadi2011boundary,yu2002natural}, let $G(\mathbf{x}, \mathbf{y})$ be  Green's function on the rectangle $\mathcal{B}$ associated with the elliptic PDE $\eqref{one_GPU:modified_helmholtz}$, which satisfies for $\mathbf{y} \in \mathcal{B}$,
\begin{align}
\triangle G(\mathbf{x}, \mathbf{y})-\kappa G(\mathbf{x}, \mathbf{y}) &= \delta(\mathbf{x}-\mathbf{y}), \quad \mathbf{x} \in \mathcal{B}, \\
G(\mathbf{x}, \mathbf{y}) &=0 \quad \mathbf{x} \in \partial \mathcal{B},
\end{align}
where $\delta(\mathbf{x} - \mathbf{y})$ is the Dirac delta function. Let $\mathbf{n}_{\mathbf{y}}$ be the unit outward normal vector at point $\mathbf{y}\in \Gamma$, and $\varphi$ be the density function. We first define the double layer boundary integral and volume integral by
\begin{align}
(W \varphi)(\mathbf{x}) := \int_{\Gamma} \frac{\partial G(\mathbf{y}, \mathbf{x})}{\partial \mathbf{n}_{\mathbf{y}}} \varphi(\mathbf{y}) d s_{\mathbf{y}}, \quad & \text { for } \mathbf{x} \in \Omega \cup \Omega^{c},\label{one_GPU:double} \\
(Yf)(\mathbf{x}) := \int_{\Omega} G(\mathbf{y}, \mathbf{x})f(\mathbf{y}) d\mathbf{y}, \quad & \text{ for } ~\mathbf{x} \in \mathbf{R}^{2}.
    \label{one_GPU:volume}
\end{align}

% Let $\psi$ be the density function, and define the single layer boundary integral by
% \begin{equation}
% (V \psi)(\mathbf{x}) := \int_{\Gamma} G(\mathbf{y}, \mathbf{x}) \psi(\mathbf{y}) d s_{\mathbf{y}} ~\text { for } ~\mathbf{x} \in \Omega \cup \Omega^{c}  \label{single}
% \end{equation}

Thanks to the symbols and properties of the involved potential and volume integral, the Dirichlet BVP $\eqref{one_GPU:modified_helmholtz}$-$\eqref{one_GPU:dirichlet_boundary}$ can be reformulated as a Fredholm boundary integral equation of the second kind\cite{kress1989linear,hsiao2008boundary} by Green’s third identity.
\begin{equation}
    \frac{1}{2}\varphi (\mathbf{x}) + (W\varphi)(\mathbf{x}) + (Yf)(\mathbf{x}) = g_{D}(\mathbf{x}), \quad \mathbf{x} \text { on } \Gamma. \label{one_GPU:second_fredholm}
    % 
    % g_D(\mathbf{x})  =(W \varphi)^{+}(\mathbf{x})+(Y f)^{+}(\mathbf{x}), \quad & \mathbf{x} \in \Gamma  \label{One_GPU:boundary_dirichlet}
\end{equation}

The solution $u(\mathbf{x})$ to the Dirichlet BVP $\eqref{one_GPU:modified_helmholtz}$-$\eqref{one_GPU:dirichlet_boundary}$ is given by 
\begin{equation}
    u(\mathbf{x}) = (W\varphi)(\mathbf{x}) + (Yf)(\mathbf{x}), \quad \mathbf{x} \in \Omega. \label{One_GPU:fredholm_dirichlet} \\
\end{equation}

Let $M>2$ be an integer and $\left\{\mathbf{x}_j\right\}_{j=0}^{\mathrm{M}}$ be a set of quasi-uniformly spaced points on the domain boundary $\Gamma$ so that each curve segment $\widetilde{\mathbf{x}_i \mathbf{x}_{i+1}}(i=0,1, \cdots, M-1)$ has nearly equal length. Numerically, the boundary integral equation $\eqref{one_GPU:second_fredholm}$ can be solved by the Richardson iteration: given an initial guess $\varphi_{0}(\mathbf{x}_{m})$, for $k \in \left\{0, 1, 2, 3, \cdots\right\}, m \in \left\{0, 1, 2, \cdots, M\right\}$, do as follows
\begin{align}
    % \hat{g}_{D}(\mathbf{x}) = g_{D}(\mathbf{x}) - (Yf)(\mathbf{x}), & \quad \mathbf{x} \in \Gamma, \\
    u_{k}^{+}(\mathbf{x}_{m}) = \frac{1}{2} \varphi_{k}(\mathbf{x}_{m}) + (W\varphi_{k})(\mathbf{x}_{m}), & \quad \mathbf{x}_{m} \in \Gamma, \label{one_GPU:richardson1} \\
    \varphi_{k+1}(\mathbf{x}_{m}) = \varphi_{k}(\mathbf{x}_{m}) + \gamma [\hat{g}_{D}(\mathbf{x}_{m}) - u_{k}^{+}(\mathbf{x}_{m})]. & \quad \mathbf{x}_{m} \in \Gamma.\label{one_GPU:richardson2}
\end{align}

Here $\hat{g}_{D}(\mathbf{x}_{m}) = g_{D}(\mathbf{x}_{m}) - (Yf)(\mathbf{x}_{m}),$ which only need to calculate once before Richardson iteration. $\mathbf{x}_{m}$ is a control node located on the boundary $\Gamma$. It can be shown that the Richardson iteration is convergent when $\gamma \in (0, 1]$. The superscript ``+'' in the BIE means one-sided limit from the domain $\Omega$. More specifically, let $w(\mathbf{x})$ be an arbitrary piecewise smooth function with discontinuities only existing at the interface $\Gamma$. We denote
 \begin{equation}
     w^{+}(\mathbf{x}) = \lim_{z \longrightarrow x, z \in \Omega} w(\mathbf{z}).
 \end{equation}
similarly, the restriction of $w(\mathbf{x})$ in $\bar{\Omega^{c}} = R^d \backslash \bar{\Omega}$, $w^{-}(\mathbf{x})$ is defined as 
\begin{equation}
    w^{-}(\mathbf{x}) = \lim_{z \longrightarrow x, z \in \bar{\Omega^{c}}} w(\mathbf{z}).
\end{equation}

Once the unknown density function $\varphi(\mathbf{x})$ is obtained when the iteration $\eqref{one_GPU:richardson2}$ converges. The unknown function $u(\mathbf{x})$ can be calculated according to the formula $\eqref{One_GPU:fredholm_dirichlet}$.
\subsection{Indirect evaluation of integrals}\label{introduce_kfbi}
In the traditional BIM method, the expression of Green's function must be explicitly known. However, the exact form of Green's function varies with the PDE, boundary condition and the domain. Although Green's function can be replaced with a neural network\cite{lin2022bigreennet} and has good numerical results in solving Laplace's and Helmholtz's equations, this method currently cannot solve the problem with variable coefficients. Within the framework of the KFBI method, there is no need to know Green's function. The integrals in $\eqref{one_GPU:double}$-$\eqref{one_GPU:volume}$ are indirectly evaluated by the equivalent interface problems. In detail, the double layer boundary integral $(W\varphi_{k})(\mathbf{x})$ and volume integral $(Yf)(\mathbf{x})$ can be written into the same form:

\begin{equation}
    \begin{array}{ll}
        \Delta v(\mathbf{x)} - \kappa v(\mathbf{x})=\mathcal{F}(\mathbf{x}), & \mathbf{x} \text { in } \Omega \cup \Omega^{c}, \\
        {[v(\mathbf{x})]=\Phi(\mathbf{x}),} & \mathbf{x} \text { on } \Gamma, \\
        {\left[\partial_{\mathbf{n}}v(\mathbf{x})\right]=0,} & \mathbf{x}\text { on } \Gamma, \\
        v(\mathbf{x})=0, & \mathbf{x} \text { on } \partial \mathcal{B}.
    \end{array} \label{one_GPU:interface}
\end{equation}

\begin{table}[ht]
    \centering
    \begin{tabular}{c|c|c}
    \hline \text { Integral } & $\mathcal{F}$ & $\Phi$  \\
    \hline
    $W\varphi$ & $\mathcal{F} = 0$ & $\Phi=\varphi$  \\
    $Yf$ & $\mathcal{F}= \tilde{f}(\mathbf{x}) = \begin{cases}f(\mathbf{x}) & \text { in } \Omega \\ 0 & \text { in } \Omega^c\end{cases}$ & $\Phi = 0$ \\
\hline
\end{tabular} \label{tab:my_label}
\end{table}
Here, [$v(\mathbf{x})$] and [$\partial_{\mathbf{n}}v(\mathbf{x})$] represent the jumps of unknown $[v(\mathbf{x})] = v^{+}(\mathbf{x}) - v^{-}(\mathbf{x})$ and its normal derivatives $[\partial_{\mathbf{n}}v(\mathbf{x})] = \partial_{\mathbf{n}}v^{+}(\mathbf{x}) - \partial_{\mathbf{n}}v^{-}(\mathbf{x}) $ respectively, $\tilde{f}(\mathbf{x})$ is the zero extension of given function $f(\mathbf{x})$. 

During the discretization of the interface problem, the discrete linear system of the interface problem $\eqref{one_GPU:interface}$ has to be corrected at the irregular nodes due to the presence of the interface $\Gamma$. The correction process needs to calculate jumps, such as [$v$], [$v_{x}$], [$v_{y}$], [$v_{xx}$], [$v_{xy}$], [$v_{yy}$] for second-order discretization, as well as modifying function values, as described in section $\ref{one_GPU:correct}$.

The jumps calculated above not only requires a correction for the discrete system in $\eqref{one_GPU:interface}$, but also interpolation of the grid-based solution $\mathbf{v}_{h}$ on the boundary. In summary, the second-order finite-difference method for solving interface is described in Algorithm $\ref{one_GPU:algorithm1}$:
\begin{algorithm}[ht]
\caption{Second-order finite difference method for interface problem $\eqref{one_GPU:interface}$}
\begin{algorithmic}[1]
\State Initialize the Cartesian grid of bounded box $\mathcal{B}$.
\State Partition the interface $\Gamma$ by quasi-uniformly control points.
\State Discretize the interface problem $\eqref{one_GPU:interface}$ with second-order finite difference method.
\State Compute jumps and correct the $\tilde{f}(\mathbf{x})$ at the irregular nodes.
\State Solve the linear system by FFT-based or geometric multigrid fast solvers.
\State Interpolate the solution to get one-sided boundary value.
\end{algorithmic} \label{one_GPU:algorithm1}
\end{algorithm}

The first step is to partition box $\mathcal{B}$ into Cartesian grid and divide the Cartesian grid nodes into regular and irregular nodes according to the boundary location. As shown in Fig.$\ref{fig:oneGPU:irregular_node}$, we define the irregular points if some of their adjacent grid nodes go cross the boundary. Black squares denote the interior irregular nodes while blue triangles denote exterior. Others are regular nodes.

For procedure of implementing steps 2-6 in Algorithm $\ref{one_GPU:algorithm1}$ on the CPU, we recommend reading reference \cite{ying2007kernel} for detail. In following section, we will focus on explaining how to execute algorithms 3-6 on single GPU in section $\ref{oneGPU}$ and multiple GPUs in section $\ref{multi_GPU}$.
% Ying provides a strategy for selecting pseudo-uniformly distributed control points in the second step$\cite{ying2013kernel}$.
% Ying provides a strategy for selecting pseudo-uniformly distributed control points in the second step$\cite{ying2013kernel}$.
\section{Single GPU algorithm} \label{oneGPU}

\subsection{Correction} \label{one_GPU:correct}
The KFBI method involves making corrections to each irregular node \cite{ying2007kernel}. Therefore, on a single GPU, each thread block is comprised of 1,024 threads, and each thread corresponds to precisely one irregular node for correction. 

\begin{figure}[htpt!]
    \centering
    \subfigure[Intersection nodes]{\includegraphics[width=0.35\textwidth]{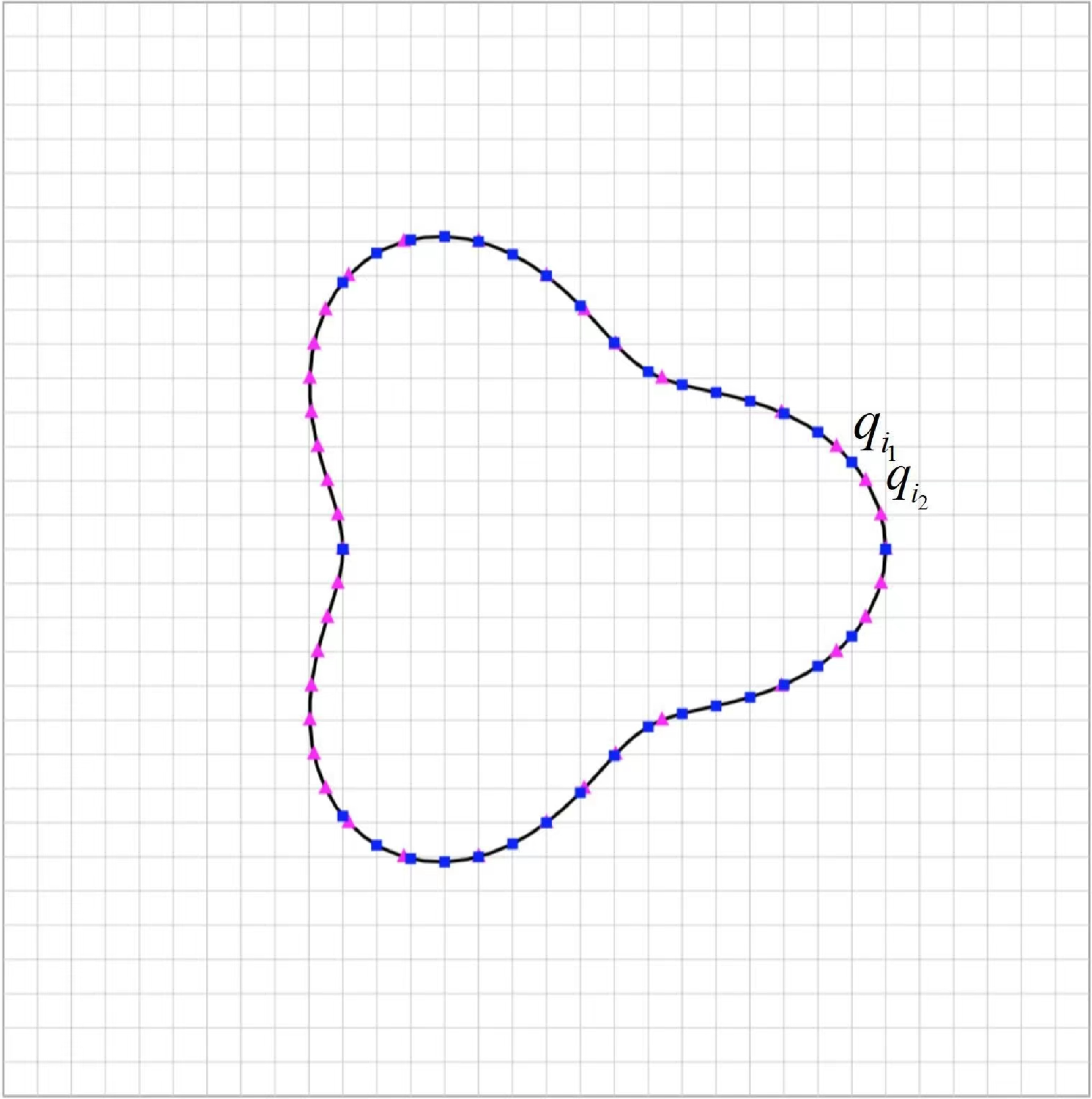}\label{fig:oneGPU:intersect_node}}
    \subfigure[Irregular nodes]{\includegraphics[width=0.55\textwidth, height=0.3501\textwidth]{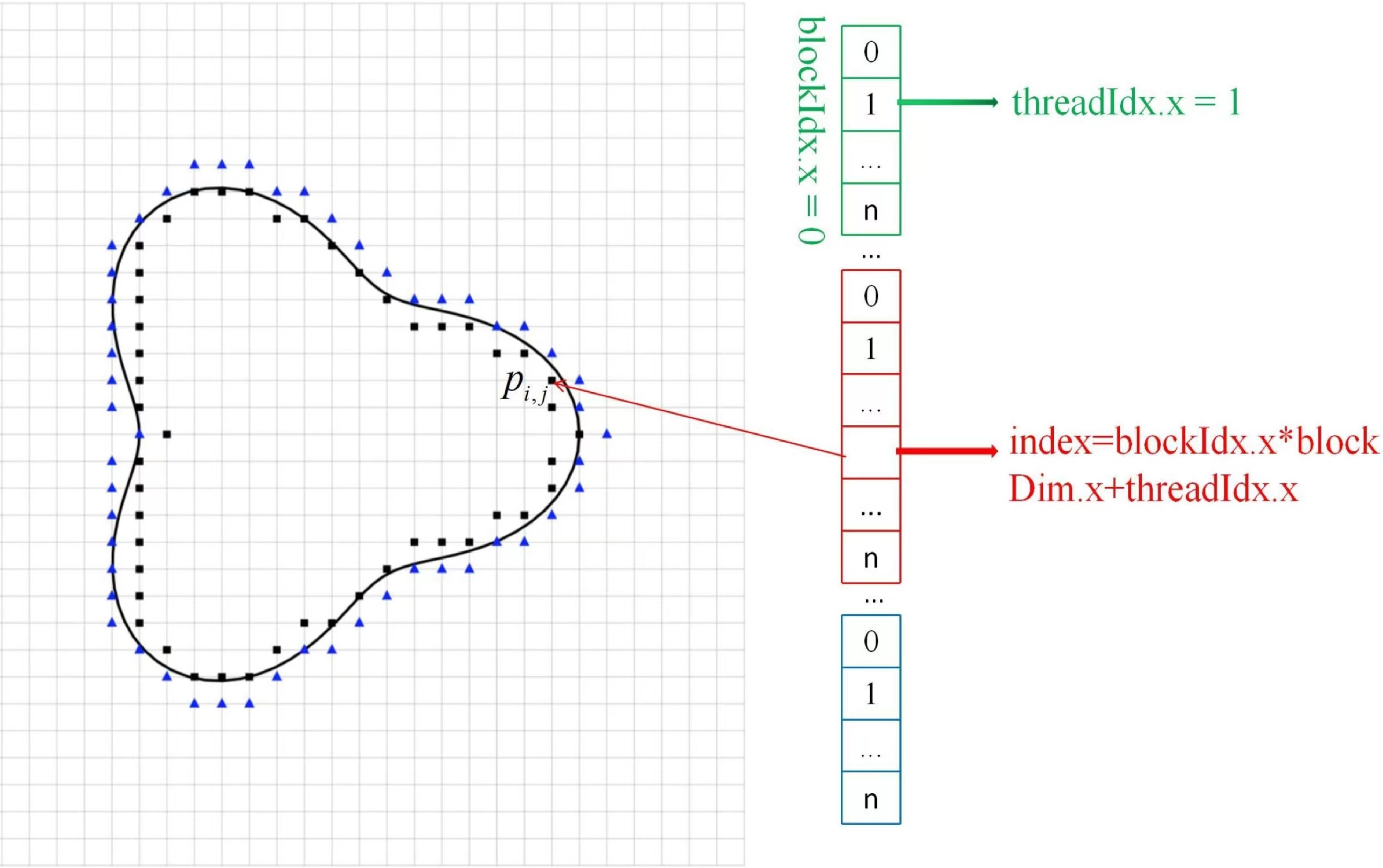}\label{fig:oneGPU:irregular_node}}
    % \subfigure[Intersection nodes]{\includegraphics[width=0.4\textwidth]{figure/one_gpu_irregular_node.eps}\label{fig:oneGPU:intersect_node}}
    % \subfigure[Irregular nodes]{\includegraphics[width=0.56\textwidth, height=0.402\textwidth]{figure/one_gpu_irregular_node_2.eps}\label{fig:oneGPU:irregular_node}}
    \caption{A graphical scheme for intersection nodes $\ref{fig:oneGPU:intersect_node}$ and irregular nodes $\ref{fig:oneGPU:irregular_node}$. In the left, pink triangles denote intersection nodes with $x$ direction while blue squares refer to intersection nodes with $y$ direction. In the right, each irregular node is controlled by one thread. Blue triangle denotes exterior irregular node while black denotes inner node.}
    \label{fig:One_GPU:irregular_node}
\end{figure}

Suppose that the rectangle domain $\mathcal{B}$ is partitioned into a uniform Cartesian grid with nodes $ \left\{\mathbf{p}_{i, j} :=\left(x_i, y_j\right): 0 \leq i \leq I, 0 \leq j \leq J\right\}$. For simplicity, assume the grid has the same spacing in the horizontal and vertical directions and denote by $h>0$ the spacing parameter, i.e., $h=x_{i+1}-x_i=y_{j+1}-y_j$. Let $v_{i,j} = v_{h}(\mathbf{p}_{i,j})$ be the finite difference approximation of $v(\mathbf{x})$ at the grid node $\mathbf{p}_{i, j}$.
One can describe the GPU version of the correction method as the following Algorithm $\ref{one_GPU:correction_algo}$ and schematic plot in Fig.\,$\ref{fig:One_GPU:irregular_node}$.
\begin{algorithm}
\renewcommand{\algorithmicrequire}{\textbf{Input:}}
\renewcommand{\algorithmicensure}{\textbf{Output:}}
\caption{Correction Procedure}
\begin{algorithmic}[1]
\Require $\text{intersection nodes, irregular nodes, } \Phi, \mathcal{F}$.
\Ensure $\text{corrected right hand side $\tilde{f}_{i,j}$ on each irregular nodes}$.
\State Locate index of irregular nodes by $\text{index} \xleftarrow[]{}$ blockIdx.x $\times$ blockDim.x + blockIdx.x
\State For irregular nodes $\mathbf{p}_{i,j}$, find the corresponding intersection nodes set $\mathcal{Q} = \left\{\mathbf{q}_{i_{1}}, \mathbf{q}_{i_{2}}, \cdots\right\}$
    \For {$\text{ each } \mathbf{q}_{i_{k}} \text{ in } \mathcal{Q}$}
    \State Interpolate $\Phi$ and calculate jumps of derivatives at of $\mathbf{q}_{i_{k}}$.
    \State  Evaluate and modify correction value $\tilde{f}_{i,j}$ by the discrete scheme on irregular nodes $\mathbf{p}_{i,j}$.
    \EndFor
\end{algorithmic}\label{one_GPU:correction_algo}
\end{algorithm}

\subsection{Interpolation} \label{one_GPU:interpolation}
% Let $v_{h}^{+}, v_{h}^{-}$ be the approximate solution of $v^{+}, v^{-}$ respectively. Assume approximate solution $v_{h}(\mathbf{x})$ is a piecewise smooth function although only the value defined at grid nodes of $\mathcal{T}_{h}$ are known. Assume the solution $v_{h}(\mathbf{x})$ is smooth in $\mathcal{B}\setminus \Gamma$. 

Fig.$\ref{interpolate}$ shows the six-point stencil of interpolated point $\mathbf{z}_{k}$ located in the shadow region of a grid cell. The stencil consists of 6 points and can be obtained by rotation or reflection transformation if $\mathbf{z}_{k}$ in another grid cell. As shown in Fig.$\ref{interpolate}$, $\mathbf{p}_{i,j}$ are the grid points for extracting value at a point $\mathbf{z}_{k} \in \Gamma$. The corresponding algorithm can be described in Algorithm $\ref{one_GPU:interpolation_algo}$:
\begin{figure}[hbpt!]
    \centering
    \includegraphics[width = 0.62\linewidth]{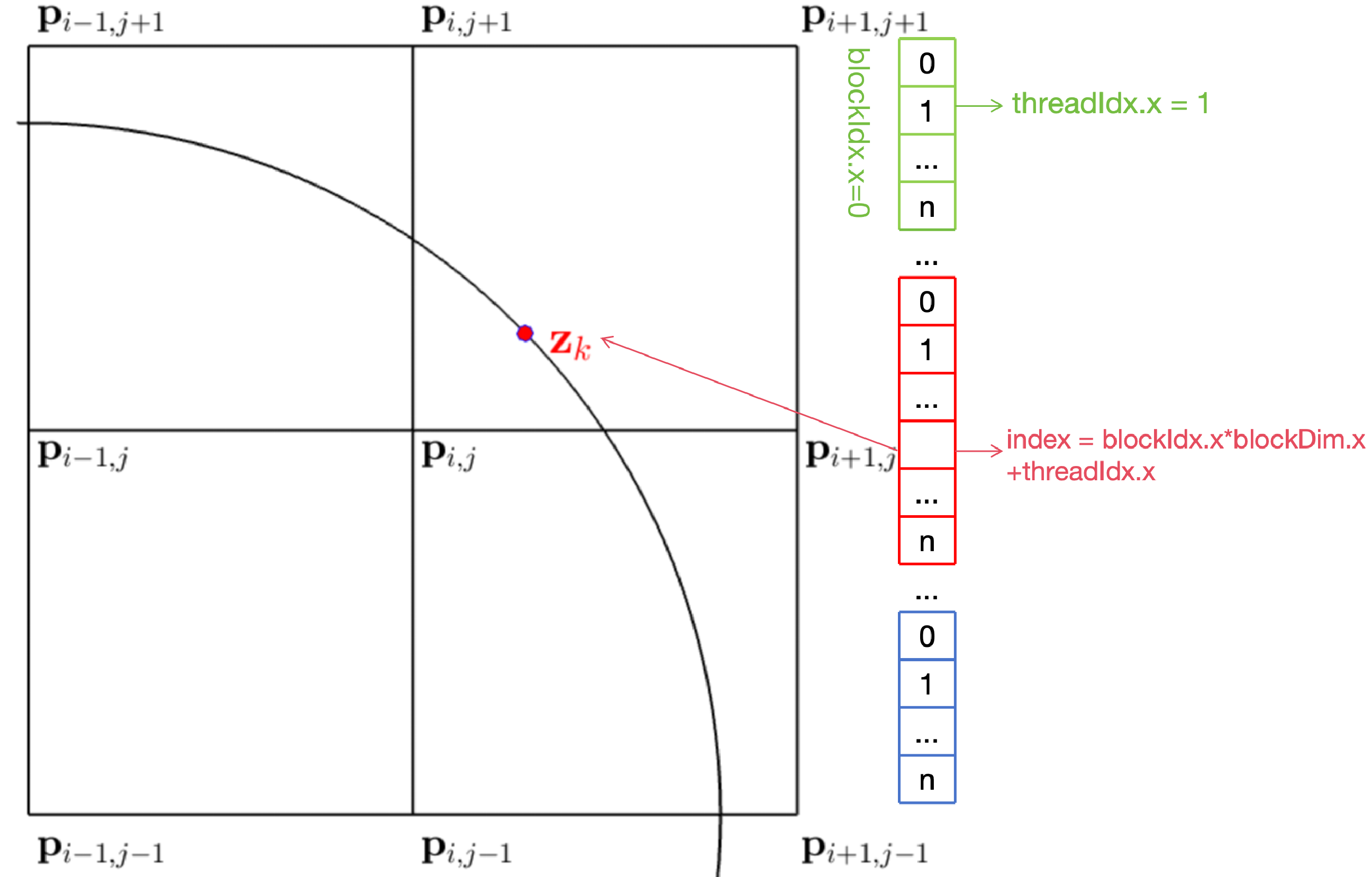}
    \caption{Interpolation schematic diagram. The evaluation of the control node $z_{k}$ depends on neighbor grid value $\mathbf{p}_{I,J}, I = \left\{i-1, i, i+1\right\}, J = \left\{j-1, j, j+1\right\}$. For parallelization, every control node is controlled by one thread. }
    \label{interpolate}
\end{figure}
\begin{algorithm}
\renewcommand{\algorithmicrequire}{\textbf{Input:}}
\renewcommand{\algorithmicensure}{\textbf{Output:}}
\caption{Interpolation Procedure}
\begin{algorithmic}[1]
\Require $\text{intersection nodes, irregular nodes, } \Phi, \mathcal{F}$.
\Ensure $\text{the value and its derivatives at control node $\mathbf{z}_{k}$}$.
\State Locate index of control nodes by $\text{index} \xleftarrow[]{}$ blockIdx.x $\times$ blockDim.x + blockIdx.x.
\For {\text{control node $\mathbf{z}_{k} \in \Gamma$}: }
    \State Compute corresponding jumps of value and its derivatives respectively.
    \State Find interpolate stencil and formulate interpolate linear system.
    \State Solve linear system and extract the boundary data on $\mathbf{z}_{k}$.
\EndFor 
\end{algorithmic}\label{one_GPU:interpolation_algo}
\end{algorithm}

\subsection{GMRES iteration with FFT-based solver} \label{one_GPU:solver} 
$\bullet \textbf{ FFT-based solver }$
% \begin{algorithm}[h!]
% \renewcommand{\algorithmicrequire}{\mathbf{Input:}}
% \renewcommand{\algorithmicensure}{\mathbf{Output:}}
% \caption{GMRES method}
% \begin{algorithmic}[1]
% \Require $\text{corrected value} \tilde{f}_{i,j}$.
% \Ensure $\text{the value of}$ $v_{i,j}$.
% \While {convergence}
%     \State \textcolor{red}{$r_{0} = \hat{g}_{D} - \mathcal{K}_{D} \varphi$}
%     \State \textcolor{blue}{$\beta = ||r_{0}||_{2}$}
%     \State \textcolor{green}{$\mu_{1} = \frac{r_{0}}{\beta}$}
%     \For{$j = 1$ to $m$} 
%         \State \textcolor{red}{$w_{j} = \mathcal{K}_{D}\mu_{j}$}
%         \For {$i = 1$ to $j$}
%             \State \textcolor{cyan}{$h_{i,j} = (w_{j}, \mu_{i}$)}
%             \State \textcolor{purple}{$w_{j} = w_{j} - h_{i,j}\mu_{i}$}
%         \EndFor
%         \State \textcolor{blue}{$h_{j+1, j} = ||w_{j}||_{2}$}
%         \State \textcolor{green}{$\mu_{j+1} = \frac{w_{j}}{h_{j+1, j}}$}
%     \EndFor
%     \State $\text{set } M_{m} = [\mu_{1}, \cdots, \mu_{m}] \text{, and } \hat{H}_{m} = (h_{i,j}) \text{ is upper Hessenberg matrix of size } (m+1) \times m$
%     \State $\text{solve a least-square problem: } \min_{y\in \mathbf{R}^{2}}||\beta \mathbf{e}_{1} - \hat{H}_{m}y||_{2}$
%     \State \textcolor{orange}{$\varphi_{m} = \varphi_{0} + M_{m}y_{m}$}
%     \If{\textcolor{blue}{$||\hat{g}_{D} - \mathcal{K}\varphi_{m}||_{2} < \epsilon$}}
%     \State $\text{convergence = true}$
%     \EndIf
%     \State $\varphi_{0} = \varphi_{m}$
% \EndWhile
% \end{algorithmic}\label{one_GPU:GMRES_algo}
% \end{algorithm}
The most important feature of the KFBI method is the conversion of a volumn or boundary integral into an interface problem. The solution to the interface problem consists of two steps. The first step is to make corrections to ensure accuracy, and the second is to solve it using a fast algorithm, such as  FFT-based solver\cite{Wu2013}. The discrete linear system  can be solved using the FFT-based solver on GPU by implementing CUDA programming and invoking cusparse\cite{naumov2010cusparse} and cufft libraries\cite{govindaraju2008high}. The algorithm can be described as 
Algorithm $\ref{one_GPU:FFT_algo}$:
\begin{algorithm}[ht]
\renewcommand{\algorithmicrequire}{\textbf{Input:}}
\renewcommand{\algorithmicensure}{\textbf{Output:}}
\caption{FFT-based solver of modified Helmholtz equation in GPU}
\begin{algorithmic}[1]
\Require $\text{corrected value}~\tilde{f}_{i,j}$.
\Ensure $\text{the solution of interface problem:} $ $v_{i,j}$.
\State Take the FFT transform on the $\tilde{f}_{i,j}$ on $y$-dimension and get transformed $\hat{\tilde{f}}_{i, j}$.
\State Solve tridiagonal linear system with right hand side term $\hat{\tilde{f}}_{i, j}$, and get solution $\hat{v}_{i, j}$. 
\State Take the inverse FFT transform on $\hat{v}_{i, j}$ on $y$-dimension, and get final solution $v_{i,j}$.
\end{algorithmic}\label{one_GPU:FFT_algo}
\end{algorithm}

In Algorithm $\ref{one_GPU:FFT_algo}$, it should be pointed out that FFT transforms depend on the boundary condition on $\partial \mathcal{B}$. If $\partial \mathcal{B}$ is the periodic, Dirichlet or Neumann, one needs to do periodic, Fast Sine Transform(FST), or Fast Cosine Transform(FCT), respectively. There is no difference in the influence of the results in the three scenarios. However, in this paper, we perform FFT on one dimension and solve the tridiagonal resulting system on another, which is the fastest and most efficient. \\
$\bullet \textbf{ GMRES iteration }$ GMRES is an iterative method for solving nonsymmetric linear systems. The method aims to approximate the solution by the vector in a Krylov subspace with minimal residual. The condition number of BIE $\eqref{one_GPU:second_fredholm}$ is relatively small, allowing for fast convergence of GMRES iteration \cite{saad2003iterative}. 

%Numerical experiment shows that compared to Richardson iteration, using GMRES to solve equation $\eqref{one_GPU:modified_helmholtz} $will require fewer iterations and therefore achieve faster convergence speed \cite{ying2013kernel, ying2014kernel, xie2019fourth}.

% Sometimes, GMRES does not converge or requires too many iterations to meet the convergence criteria. Therefore, in most cases, GMRES must include a preconditioning step to improve its convergence. 

From $\eqref{one_GPU:richardson1}$, we denote  : 
\begin{align}
    \mathcal{K}_{D}(\varphi)(\mathbf{x}) & := (\frac{1}{2}I + W)(\varphi)(\mathbf{x}), \quad \mathbf{x} \in \Gamma
\end{align}
\begin{algorithm}[ht]
\renewcommand{\algorithmicrequire}{\textbf{Input:}}
\renewcommand{\algorithmicensure}{\textbf{Output:}}
\caption{GMRES with restarts in GPU}
\begin{algorithmic}[1]
\Require $\text{corrected value} ~\tilde{f}_{i,j}$.
\Ensure $\text{the solution of BIE}($\ref{one_GPU:second_fredholm}): $\varphi_{m}$.
\State convergence = false
\While {convergence == false}
    \State $r_{0} = \hat{g}_{D} - \mathcal{K}_{D} \varphi$
    \State $\beta = ||r_{0}||_{2}$
    \State $\mu_{1} = \frac{r_{0}}{\beta}$
    \For{$j = 1$ to $m$} 
        \State $w_{j} = \mathcal{K}_{D}\mu_{j}$
        \For {$i = 1$ to $j$}
            \State $h_{i,j} = (w_{j}, \mu_{i}$)
            \State $w_{j} = w_{j} - h_{i,j}\mu_{i}$
        \EndFor
        \State $h_{j+1, j} = ||w_{j}||_{2}$
        \State $\mu_{j+1} = \frac{w_{j}}{h_{j+1, j}}$
    \EndFor
    \State $\text{Set } M_{m} = [\mu_{1}, \cdots, \mu_{m}] \text{, and } \hat{H}_{m} = (h_{i,j}) \text{ is upper Hessenberg matrix of order } (m+1) \times m$
    \State $\text{Solve a least-square problem: } \min_{y\in \mathbf{R}^{2}}||\beta \mathbf{e}_{1} - \hat{H}_{m}y||_{2}$
    \State $\varphi_{m} = \varphi_{0} + M_{m}y_{m}$
    \If{$||\hat{g}_{D} - \mathcal{K}\varphi_{m}||_{2} < \epsilon$}
    \State $\text{convergence = true}$
    \EndIf
    \State $\varphi_{0} = \varphi_{m}$
\EndWhile
\end{algorithmic}\label{one_GPU:GMRES_algo}
\end{algorithm}
The main points of the GMRES method with restarts are presented in Algorithm $\ref{one_GPU:GMRES_algo}$. It is worth noting that the matrix-vector product in GMRES is replaced by the KFBI method(line 2 and 6). As for computing vector norm(line 3, 11, and 17), inner product(line 8), scalar-vector operation(line 4 and 12), and some axpy operation(line 9), they are implemented by the cuda kernel functions. Since the calculated $\varphi$ need to be reintegrated into the next iteration, the GPU information needs to be transferred back to the CPU after each GMRES iteration to obtain updated $M_{m-1}$ to $M_{m}$(line 14). Therefore, the calculations for the least squares method are also performed on the CPU(line 15).
\section{Multi-GPUs algorithm} \label{multi_GPU}
\subsection{Domain decomposition}

While a single GPU has demonstrated commendable performance in numerous applications, the computational demands entailed in simulating extensive more large-scale problems surpass the capabilities of a single GPU. It is necessary to explore the development of parallelization techniques using multiple GPUs to address this issue. For descriptive convenience, our multiple-GPU algorithm is presented using a two-dimensional grid as a paradigm, with the situation for three-dimensional grids being analogous.

Our study employs the domain decomposition method to partition a 2D grid with $Nx \times Ny$ dimensions into $m$ parts along the $x-$coordinate direction. Each part represents a sub-domain; the total number of sub-domains is denoted as $m, m$ equals the total number of GPUs in use. Each sub-domain is assigned to a corresponding process equipped with a GPU for computation. Throughout the simulation, the variables of each sub-domain persistently reside in the global memory of the assigned GPU.
\begin{figure}[ht]
    \centering{\includegraphics[width=0.4\textwidth]{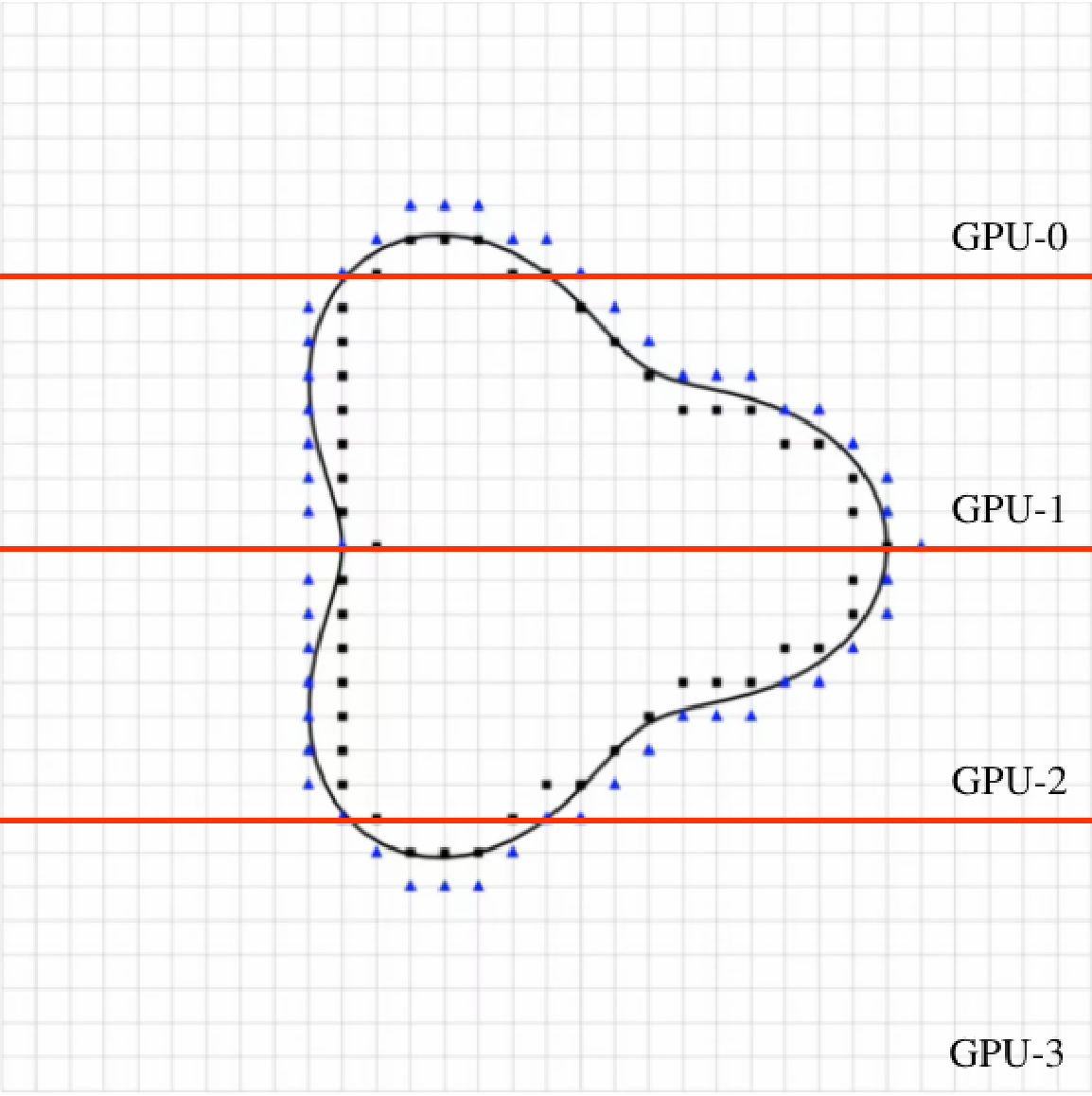}}
    \caption{Partition the Cartesian grid into four subdomains along the y-axis, with each subdomain assigned to a dedicated process equipped with a GPU for computation.}
\end{figure}

In the interface problem-solving process, each process stores the relevant information of curve points within its allocated region. Initially, all the data points on the boundary are gathered to calculate the solutions at these specific curve points. Subsequently, the computed results are distributed back to the respective processes.

Before set up the iteration in Section \ref{one_GPU:solver}, the CPU's grid data and boundary data are distributed to their respective processes based on regions. Subsequently, this data is efficiently copied to the device. Throughout the iterative solving process, each process is responsible for processing the pertinent information related to grid points, control points, and intersection nodes within its designated region, and all computations are carried out on the GPU. Upon completion of the iteration, the processed data is then returned to the host system. This approach ensures a standardized and optimized procedure for utilizing CUDA-enabled devices to expedite the solving process.

\subsection{Data communication}
\begin{figure}[ht]
    \centering{\includegraphics[width=0.4\textwidth]{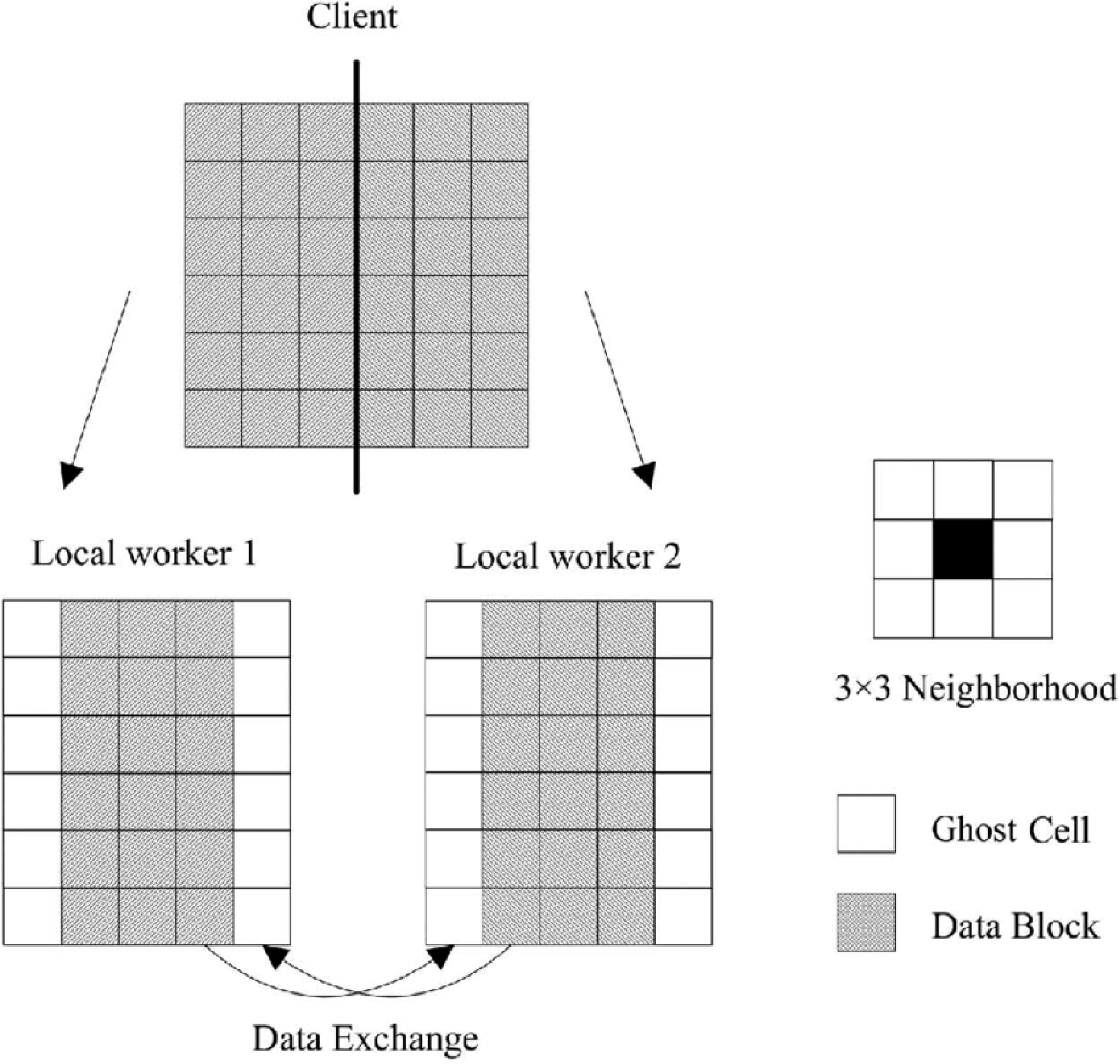}}
    \caption{Each process is accompanied by a peripheral layer of a virtual subdomain, situated at the boundaries of their respective regions, which serves the purpose of receiving boundary data transmitted by neighboring processes.}
 	 	%	\label{125000}
\end{figure}

In order to facilitate efficient data exchange, MPI is utilized on the CPU to transfer data from two layers of internal units adjacent to subdomain boundaries. We allocate memory spaces on both the host and the device to store secondary boundary data. Once the boundary of a specific subdomain is computed, each device uploads the necessary boundary data to the virtual region, which is then downloaded by neighbouring devices to prepare for the subsequent computational steps. This approach optimizes data exchange while minimizing data redundancy between the CPU and GPU.

Given that there are fewer control points than grid points, the time and computational cost associated with collecting and disseminating potential values, denoted as $\phi$ and $\psi$, across various boundary regions is relatively low. During the process of solving the interface problem \eqref{one_GPU:interface}, each process is tasked with obtaining potential function information for all control points. To streamline this procedure, MPI is employed to consolidate the potential function information before the interface problem calculation begins. Upon completion of the calculation, MPI is once again used to distribute the results back to their respective processes. 
\subsection{Poisson solver}\label{distribute Poisson}
% \subsubsection{Algorithm}
% For the Poisson equation, we employ a three-step process to solve the modified linear system:
% \begin{enumerate}
% \item  FFT in $x$-direction: Firstly, the data is transformed by FFT. This operation is performed in the $y$ direction, and the Fourier coefficients are computed using the CUFFT library. This step allows for efficient manipulation of the data in the frequency domain.
% \item Distributed Arrowhead Decomposition method: Next, we utilize a distributed arrowhead decomposition method to decompose the tridiagonal linear system in the $x$ direction. This technique facilitates solving the linear tridiagonal equations for each sub-domain on the corresponding GPU. We leverage the cusparseDgtsv2StridedBatch API from the CUSPARSE library for an efficient and accurate solution of the linear equations.
% \item IFFT in $x$-direction: Finally, the data is transformed back to the spatial domain using the inverse fast Fourier transform. This operation is performed in the $y$ direction, and the Fourier coefficients are calculated once again using the CUFFT library. This step allows us to obtain the final solution of the Poisson equation.
% \end{enumerate}
For the FFT-based solver of the Poisson equation, we still follow the process Algorithm \ref{one_GPU:FFT_algo}. Unlike single GPU, we solve tridiagonal linear equations using the distributed arrowhead decomposition method(ADM)\cite{arrowhead2017}. The ADM is an efficient algorithm, and the resulting system is suited for designing distributed algorithms for each sub-domain on the corresponding GPU.

\subsubsection{Arrowhead decomposition method}
The linear equation system to be solved is denoted as : $Au = f$.
\begin{figure}[htb]
    \centering
    \subfigure[Intial linear system]{
	\includegraphics[width=0.3\textwidth]{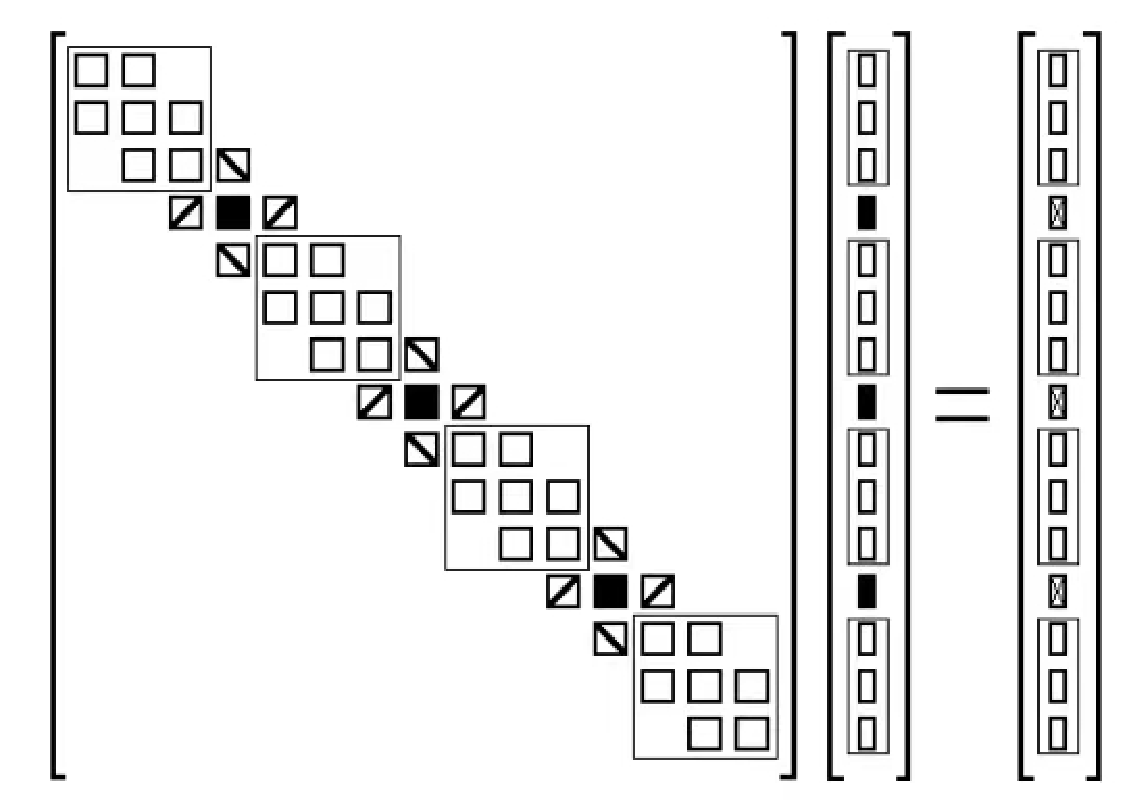}
	\label{fig:multi_GPU:triangle_a}
    }
    \subfigure[Rearranged linear system]{
   	\includegraphics[width=0.3\textwidth]{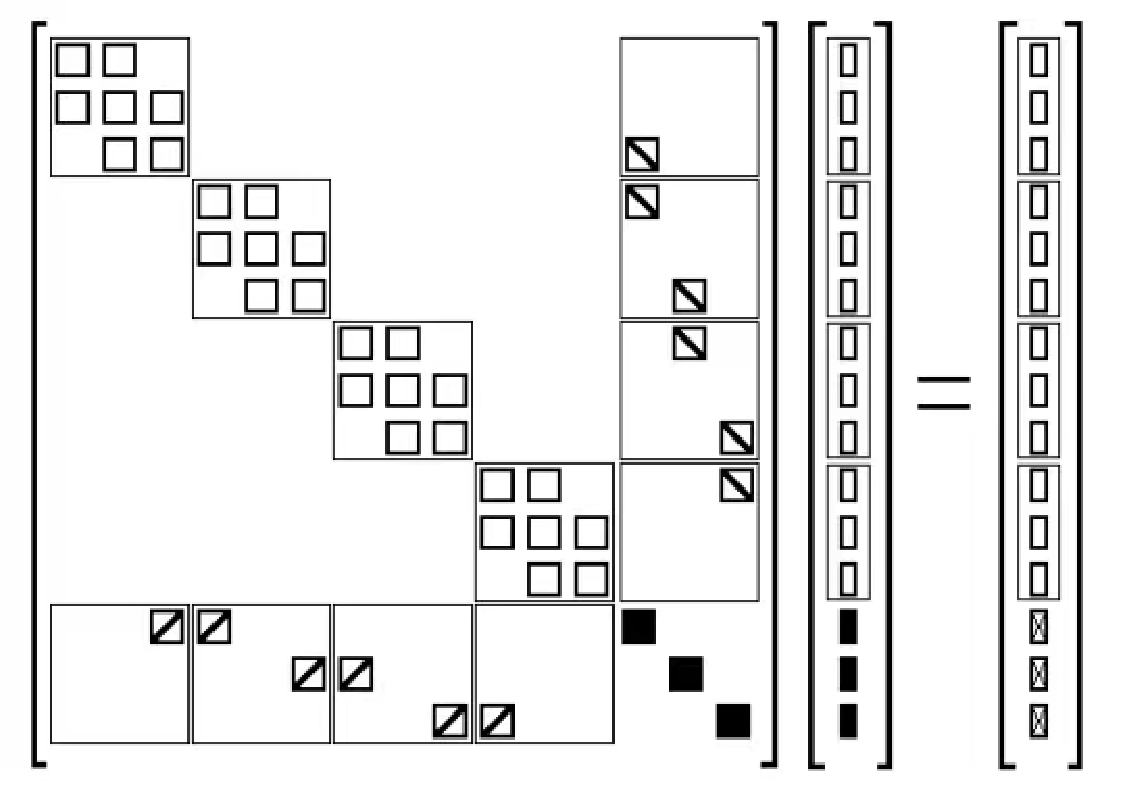}
	\label{fig:multi_GPU:triangle_b}
    }
    \subfigure[Notation of the rearranged system]{
   	\includegraphics[width=0.3\textwidth]{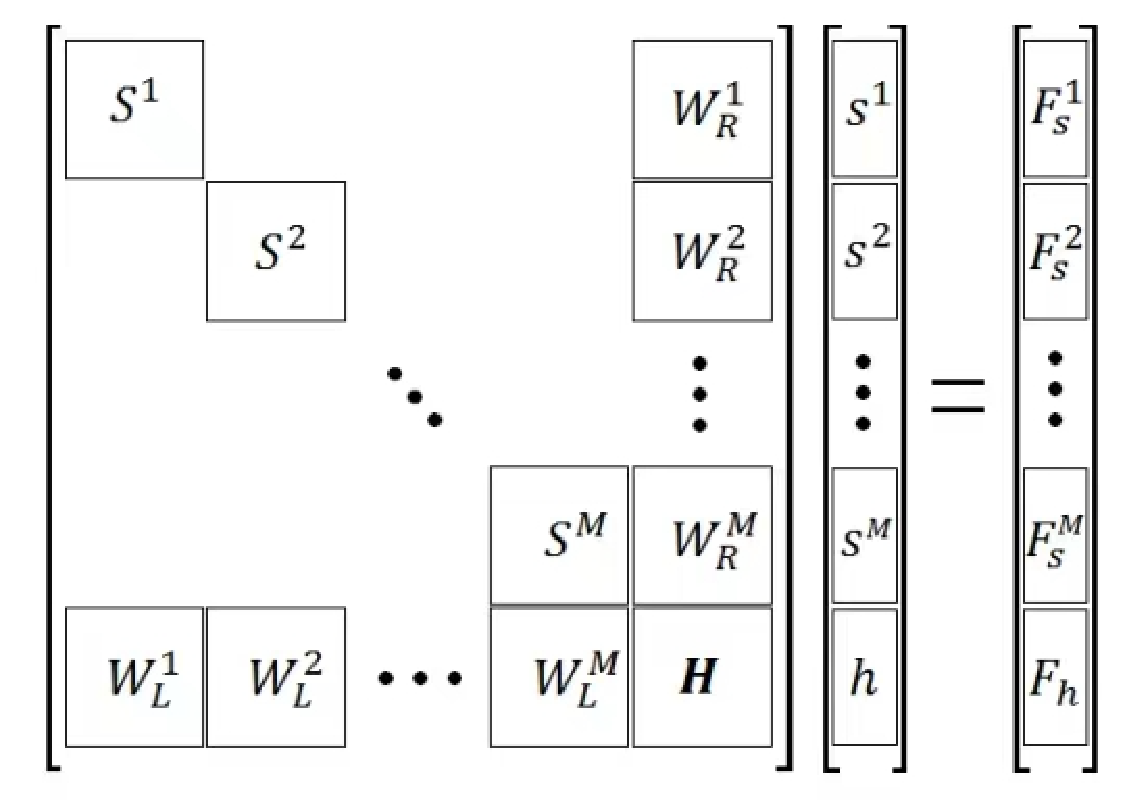}
	\label{fig:multi_GPU:triangle_c}
    }
    \caption{A graphical scheme for rearrangement of the initial block-tridiagonal linear system into the equivalent form, coming from \cite{arrowhead2017}. Left: the initial matrix $A$ with blocks. Center: after rearrangement, the initial matrix becomes "arrowhead" matrix. Right: the denotation of the "arrowhead" matrix.}
    \label{Schur Complement Method}
\end{figure}
Fig.\ref{Schur Complement Method} depicts the concept of ADM. A similarity transformation transforms the initial block-tridiagonal linear system\eqref{block-tridiagonal} into an equivalent block matrix form. This reorganization is particularly advantageous for region decomposition and the design of distributed algorithms, as each block matrix's linear system exhibits a degree of independence. The reordering is carried out by exchanging block rows and block columns, which, in turn, affects the elements of the unknown vector and the right-hand side vector. The resulting matrix structure can be represented as a $2 \times 2$ block matrix.
\begin{equation}
\left(\begin{array}{cc}
\mathbf{S} & \mathbf{W}_R \\
\mathbf{W}_L & \mathbf{H}
\end{array}\right)\left(\begin{array}{l}
s \\
h
\end{array}\right)=\left(\begin{array}{l}
F_s \\
F_h
\end{array}\right)
\label{block-tridiagonal}
\end{equation}

In this context, the unknown solution vector $h$ corresponds to the movable component of the complete solution, as illustrated in Figure \ref{Schur Complement Method}. The square super-block $\mathbf{S}$ comprises new independent tridiagonal blocks $S^k$, with $k=1, \ldots, M$, positioned along the diagonal. The matrix elements $\mathbf{H}$ form the lower-right coupled super-block, representing the coupling of unknowns at the interface. Additionally, the other horizontal super-blocks $\mathbf{W}_R$ and $\mathbf{W}_L$ are supplementary components within the matrix, signifying the internal unknowns of the coupling processors at the interface. The following relationships determine the solution of the system\eqref{block-tridiagonal}.

\begin{equation}
\left\{\begin{array}{l}
s=\mathbf{S}^{-1} F_s-\mathbf{S}^{-1} \mathbf{W}_R h \\
h=\left(\mathbf{H}-\mathbf{W}_L \mathbf{S}^{-1} \mathbf{W}_R\right)^{-1}\left(F_h-\mathbf{W}_L \mathbf{S}^{-1} F_s\right)
\end{array}\right.
\label{block-tridiagonal:schur}
\end{equation}

These relationships involve matrix products and inversions, which can be parallelized to a certain extent. The independence of blocks $S^k$ allows for distributed parallel computation of the products $\mathbf{S}^{-1} F_s$. In practice, rather than computing inversions, we efficiently solve the distributed linear systems $Sx = F_s$ due to the special properties of $S^k$. The sparse structure of $\mathbf{W}_L$ and $\mathbf{W}_R$ significantly reduces the number of matrix operations in\eqref{block-tridiagonal:schur}. Once a portion of the total solution $\mathbf{X}=(s, h)^T$ is obtained from the second relationship in\eqref{block-tridiagonal:schur}, i.e., $h=\left(h^1, \ldots, h^{M-1}\right)^T$, the remaining parts can be computed in parallel over GPU $k$ using the provided formula.
\begin{equation}
s^k=z^k-Z^k h^{k-1}-Z^k h^k
\label{trans}
\end{equation}
where formally $z^k=\left(S^k\right)^{-1} F_s, Z^k=\left(S^k\right)^{-1} W_R^k$, and $h^0=h^M=0$.

% \begin{figure}[H]
%  	 		\centering
%  	 		{
%  	 			\includegraphics[width=0.7\textwidth]{figure/trianglematrix.png}
%  	 		}
%  	 		\centering
%     \caption{A graphical scheme for rearrangement of the initial block-tridiagonal linear system
% into the equivalent form.}
%     \label{Schur Complement Method}
%  	 	%	\label{125000}
% \end{figure}

% \begin{figure}[htb]
%     \centering
%     \subfigure[Intial linear system]{
% 	\includegraphics[width=0.3\textwidth]{figure/triangleM_a.eps}
% 	\label{fig:multi_GPU:triangle_a}
%     }
%     \subfigure[Rearranged linear system]{
%    	\includegraphics[width=0.3\textwidth]{figure/triangleM_b.eps}
% 	\label{fig:multi_GPU:triangle_b}
%     }
%     \subfigure[Notation of the rearranged system]{
%    	\includegraphics[width=0.3\textwidth]{figure/triangleM_c.eps}
% 	\label{fig:multi_GPU:triangle_c}
%     }
%     \caption{A graphical scheme for rearrangement of the initial block-tridiagonal linear system into the equivalent form, coming from \cite{arrowhead2017}. Left: the initial matrix $A$ with blocks. Center: after rearrangement, the initial matrix becomes "arrowhead" matrix. Right: the denotation of the "arrowhead" matrix.}
%     \label{Schur Complement Method}
% \end{figure}

The distributed algorithm for solving tridiagonal matrices can be outlined as follows:
\begin{algorithm}[H]
\renewcommand{\algorithmicrequire}{\textbf{Input:}}
\renewcommand{\algorithmicensure}{\textbf{Output:}}
\caption{Distributed solver for arrowhead decomposition method}
\begin{algorithmic}[1]
\Require $\text{corrected value} ~\tilde{f}_{i,j}$ and known matrix $A$.
\Ensure $\text{the solution of}$ $Au = f$.
\State Decompose the coefficient matrix $A$ and the right-hand side vector $f$ into $m$ subregions, where $m$ equals the number of processes.
\State Precompute the Schur complement matrices $(H-W_LS^{-1}W_R)^{-1}, ~S^{-1}W_R$ independently. 
\State Compute $S_{k}^{-1}F_{s}, ~k = \left\{1,2, \cdots, m\right\}$ and exchange data by passing the first row of $x_k$ from the $(k+1)^{th}$ processes to the auxiliary boundary of the $k^{th}$ process to compute $\left(F_h-\mathbf{W}_L \mathbf{S}^{-1} F_s\right)$.
\State Compute $h^{k}, k = \left\{1, \cdots, m-1\right\}$ by $\eqref{block-tridiagonal}$ and $h^{0}, h^{m}$ are set to 0.
\State Evaluate $s^{k}, k = \left\{1, \cdots, m\right\}$ by formula $\eqref{trans}$, where the vector $h^{k-1}$ is passed from $(k-1)^{th}$ process  to the auxiliary boundary of the process $k$.
\end{algorithmic}
\end{algorithm}
 
 In step 1, each process is assigned the task of handling storage and computations for variables within its respective subregion. This approach ensures a standardized and optimized procedure for distributing the workload across multiple processes. In step 2, computing these matrices avoid redundant calculation during the iterations. For the vector update operations in both step 1 and 2, each vector is divided into some segments according to the number of devices. Each segment pair forms a subtask in a device and these subtasks are computed simultaneously. For the dot product in point, the vectors $\vec{x}$ and $\vec{y}$ are cut into segments and computes on the devices in parallel firstly, then the local sum was calculated using the API $reduce$ in the thrust library. The MPI is used to solve the final global sum, the vector norm can easily be calculated when the dot product is acquired. In our program, an original vector is partitioned sequentially. Any vector is stored as a one-dimensional array. A set of vectors are managed through a two-dimensional array. For the matrix-vectors product, here is actually the solution to the interface problem.

\subsection{Algorthm summary}\label{Sec:Algorthm summary}
The structure of the GPU-accelerated distributed KFBI algorithm can be found hereinafter. The
individual steps are interleaved by communication calls, as
visualized by printing the communication in italic.

\begin{itemize}
    \item Procedure 1: Initialize the Cartesian grid 
    \begin{enumerate}
        \item Use quasi-uniformly spaced points $\textbf{z}_{i}$ to discretize the interface $\Gamma$;
        \item Partition $\mathcal{B}$ into a uniform Cartesian grid $\mathcal{T}_{h}$;
        \item Identify the regular and irregular nodes of the Cartesian grid;
        \item Find intersecting points located between $\Gamma$ and Cartesian grid line.
    \end{enumerate}
    \item Procedure 2: Evaluate boundary data on control nodes
    \begin{enumerate}   
        \item The boundary point data is scattered to the corresponding process according to the region;
        \item Compute jumps of partial derivatives at control nodes;
        \item Solve interface problems by section $\ref{distribute Poisson}$;
        \item Exchange grid data between adjacent processes;
        \item Extract boundary data $u^{+}(\textbf{x})$ or $\partial_{\textbf{n}} u(\textbf{x})$ for Dirichlet or Neumann BVP respectively;
        \item \text{Collect boundary point data and calculate errors}.
    \end{enumerate}
    \item Procedure 3: The GMRES iteration
    \begin{enumerate}
        \item Choose an initial guess $\varphi_0$ or $\psi_0$ and distribute it to different GPUs to start the GMRES iteration for the Dirichlet BVP or Neumann BVP respectively;
        \item Perform the GMRES iteration according to \ref{one_GPU:solver};
        \item Repeat the previous steps 2 until the discrete residual of the boundary integral equation is sufficiently small in some norm.
    \end{enumerate}
    \label{time-dependence solver}
\end{itemize}

\section{Numerical Results} \label{result}

To study the numerical accuracy and efficiency of the methods above, in this section, we present the numerical results for the Laplace equations, the reaction-diffusion equations, and the Stokes equations in an irregular domain. The bounding box $\mathcal{B}$ embedding the domain $\Omega$ for solving the interface problem is specified as a square(cube), whose size as well as the curve(surface) parameters are given respectively in the description of each numerical example.

The following examples give the convergent error of the numerical discretization scheme. Taking two dimensions as an example, the error is defined as $e_{i j}$ with $e_{i j}=\left(u_h\right)_{i j}-\left(u^*\right)_{i j}$, where $N$ is the number of interior grid nodes, $u^*$ is the exact solution, $u_h$ is the numerical solution with step size $h$. Denote by $\left\|\mathrm{e}_h\right\|_{\infty}$ and $\left\|\mathrm{e}_h\right\|_2$ the discrete maximum norm and the scaled discrete $l_2$-norm of $e_{i j}$ respectively, i.e.,
$$
\begin{aligned}
\left\|\mathbf{e}_h\right\|_{\infty} & =\max _{\left(x_i, y_j\right) \in \Omega}\left|e_{i j}\right| \\
\left\|\mathbf{e}_h\right\|_2 & =\sqrt{\frac{1}{N} \sum_{\left(x_i, y_j\right) \in \Omega}\left|e_{i j}\right|^2}
\end{aligned}
$$

To check the algorithm's accuracy, we verify the numerical error in each case with the grid refinement. The GMRES iteration stops when the iterated residual in the discrete $\ell^2$-norm relative to that of the initial residual is less than a prescribed tolerance and is fixed to be $10^{-8}$. The corresponding table for each case lists the step size, the number of grid points, the CPU times, the GPU times, and the speedup ratios. Numerical results on the Cartesian grid to the problem are also displayed in the plots for each fixed time.

In addition, we perform numerical experiments on eight NVIDIA GeForce RTX 3090 graphics cards, which contain 10496 cores organized in 84 streaming multiprocessors (MPs). Moreover, it provides 24GB of device memory with a memory bandwidth of 936GB/s, accessible by all its cores and the CPU through the Intel(R) Xeon(R) Gold 6330 CPU with 28 cores.

\begin{table}[H]
\centering
\begin{tabular}{c|c|ccccc}
\hline Poisson & \text { grid size }  & $512^2$ & $1024^2$ & $2048^2$ &  $4096^2$\\
\hline FFT
 & GPU time & $0.19 \mathrm{~s}$ & $0.27 \mathrm{~s}$ & $0.66 \mathrm{~s}$ & $1.65\mathrm{~s}$  \\
\hline Multigrid&  GPU time & $0.25 \mathrm{~s}$ & $0.61 \mathrm{~s}$ & $2.08 \mathrm{~s}$ & $8.66\mathrm{~s}$  \\
\hline
\end{tabular}
    \caption{\footnotesize{\small{Comparison of different Poisson solvers on DBVP of the Laplace equation.}}}
    \label{tab:FFTvsMultigrid}
\end{table}
\begin{table}[H]
\centering
\begin{tabular}{c|c|ccccc}
\hline Iteration & \text { grid size }  & $512^2$ & $1024^2$ & $2048^2$ &  $4096^2$\\
\hline Richardson& \text {CPU time}& $5.62 \mathrm{~s}$ & $23.11 \mathrm{~s}$ & $92.05 \mathrm{~s}$ & $380.12\mathrm{~s}$ \\
 & GPU time & $0.19 \mathrm{~s}$ & $0.27 \mathrm{~s}$ & $0.66 \mathrm{~s}$ & $1.65\mathrm{~s}$  \\
\hline GMRES&  \text {CPU time}& $1.38 \mathrm{~s}$ & $5.81 \mathrm{~s}$ & $25.92 \mathrm{~s}$ & $110.21\mathrm{~s}$ \\
 & GPU time & $0.13 \mathrm{~s}$ & $0.16 \mathrm{~s}$ & $0.18 \mathrm{~s}$ & $0.25\mathrm{~s}$  \\
 \hline BiCGSTAB&  \text {CPU time}& $1.99 \mathrm{~s}$ & $8.01 \mathrm{~s}$ & $34.47 \mathrm{~s}$ & $147.49\mathrm{~s}$ \\
 & GPU time & $0.45 \mathrm{~s}$ & $0.45 \mathrm{~s}$ & $0.63 \mathrm{~s}$ & $1.01\mathrm{~s}$  \\
\hline
\end{tabular}
    \caption{\footnotesize{Comparison of various iterative methods for solving the Dirichlet boundary value problem (DBVP) associated with the Laplace equation, with a fixed tolerance level of $1e-08$ for the Richardson, GMRES, and BiCGSTAB methods.}}
    \label{tab:iter}
\end{table}
\subsection{Single GPU results}
\textbf{Example 1.} The results presented in Tab.\,\ref{tab:FFTvsMultigrid} clearly demonstrate that, within the specified range of simulation scales, the FFT+tridiagonal Poisson solver outperforms the geometric multigrid solver in terms of efficiency. This observation is consistent with the findings reported in\cite{Gholami2016}. Analyzing the results from Tab.\,\ref{tab:iter}, it is evident that the GMRES method achieves the lowest iteration number, resulting in reduced CPU and GPU time consumption. As a result, for subsequent numerical experiments, this study adopts the parallel FFT+tridiagonal Poisson solver in combination with the GMRES method for computational purposes.

This example solves the boundary value problem of the Laplace equation on the circle domain(the parameters $r_a=1.0$, $r_b=1.0$) and the rotated star-shaped domain(the parameters $m=4.0,6.0,8.0$, $r=1.0$, $c=0.2$), with the Dirichlet boundary condition. The boundary conditions are chosen so that the exact solution reads
$$
\begin{gathered}
u(x, y)=\exp(x)\cos(y) + \exp(y)\sin(x)
\end{gathered}
$$

\begin{figure}[htpt!]
    \centering{\includegraphics[width=0.6\textwidth]{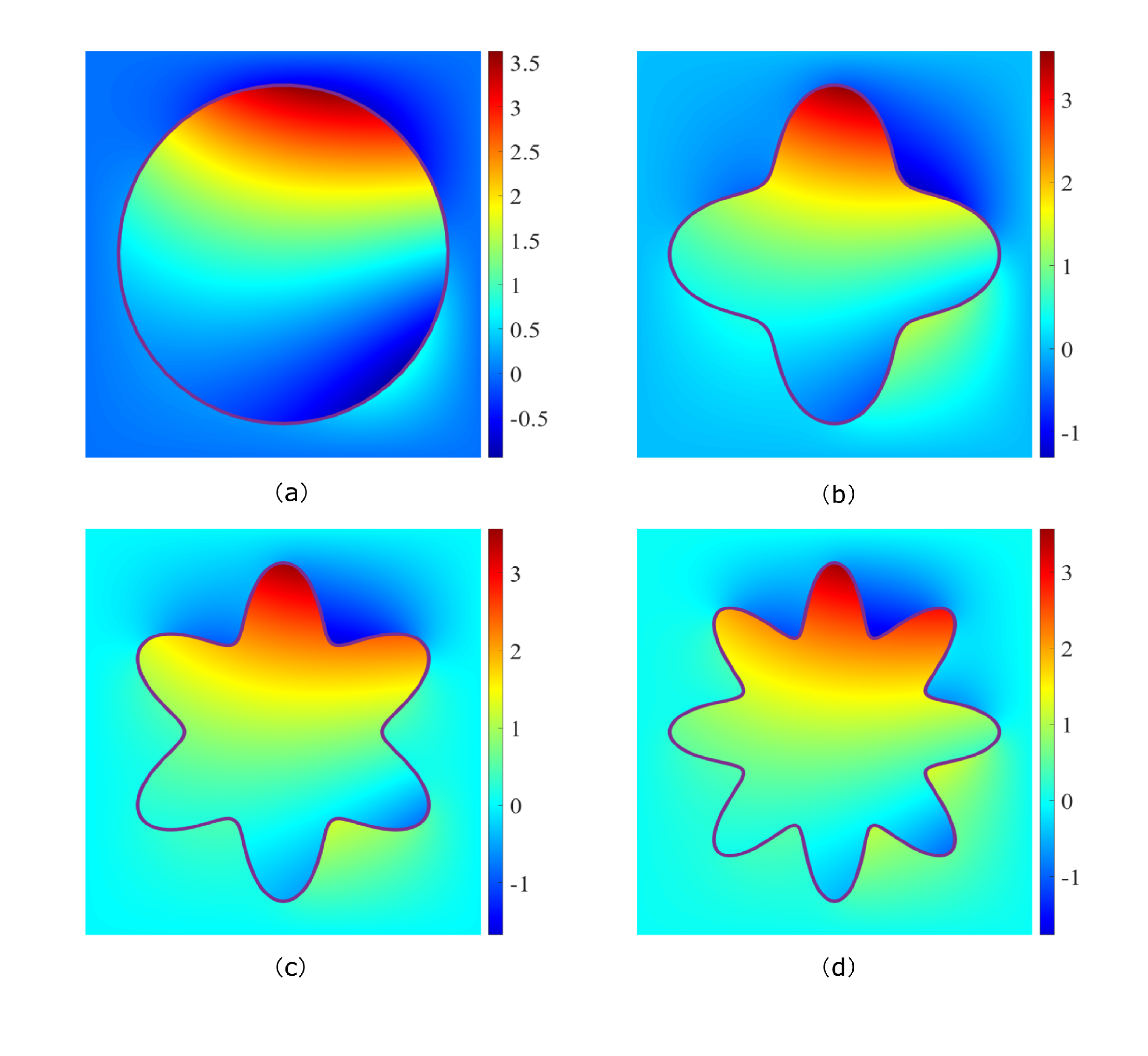}}
    \caption{The numerical solutions on the circle and star-shaped domain on the $2048 \times 2048$ grid.(a)Circle domain. The radius r = 1.0. (b)The star-shaped domain. The fold number m = 4.0, radius r = 1.0, c = 0.2. (c)The star-shaped domain. The fold number m = 6.0, radius r = 1.0, c = 0.2. (d)The star-shaped domain. The fold number m = 8.0, radius r = 1.0, c = 0.2.}
    \label{One_Poisson1}
\end{figure}

The bounding box $\mathcal{B}$ for the interface problem is set to be $\mathcal{B}=(-1.2,1.2) \times(-1.2,1.2)$. Numerical results are plotted in Fig.\,\ref{One_Poisson1}.

\begin{figure}[htb]
    \centering
    \subfigure[The velocity field $u^{(1)}$]{\includegraphics[width=0.3\textwidth]{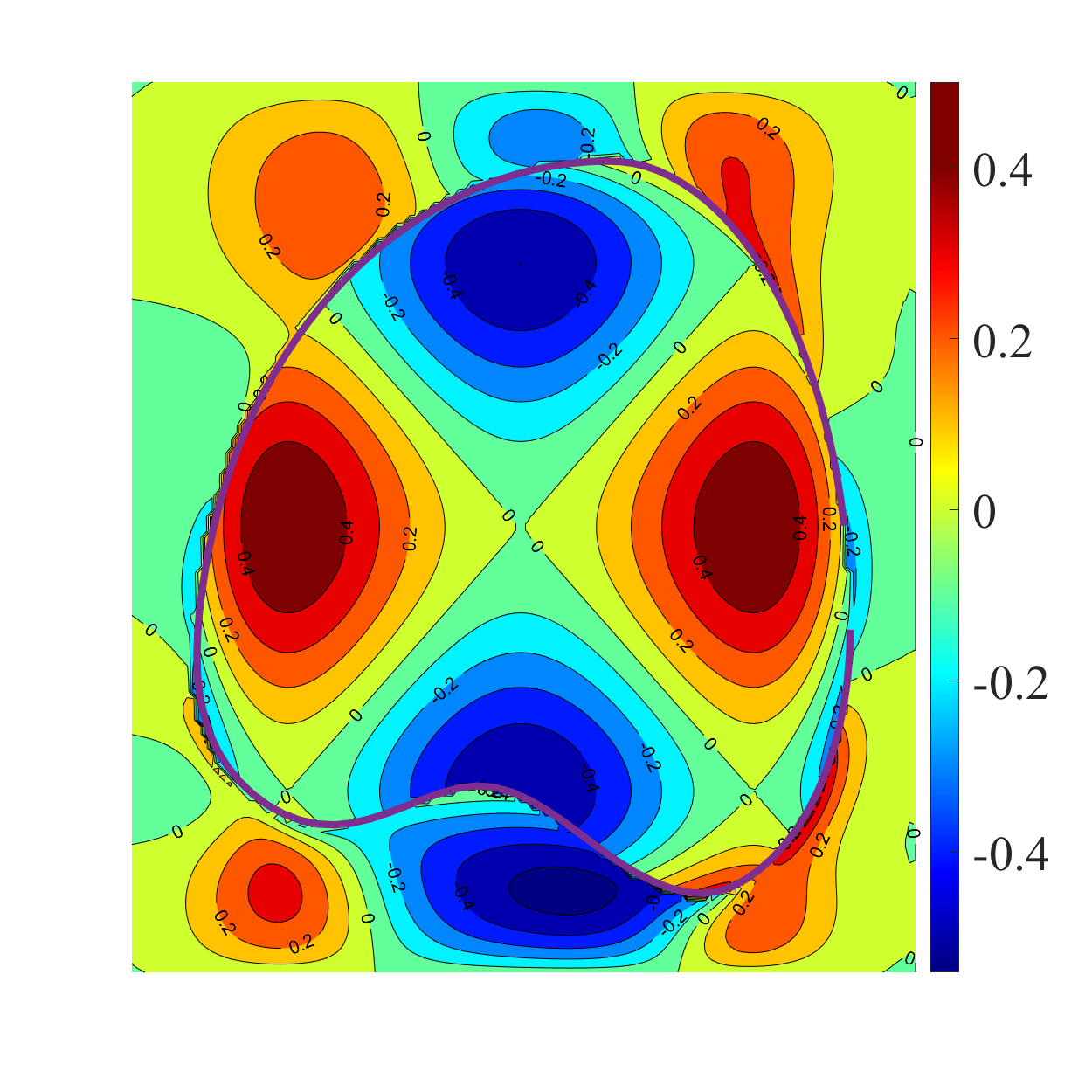}}
    \subfigure[The velocity field $u^{(2)}$]{\includegraphics[width=0.3\textwidth]{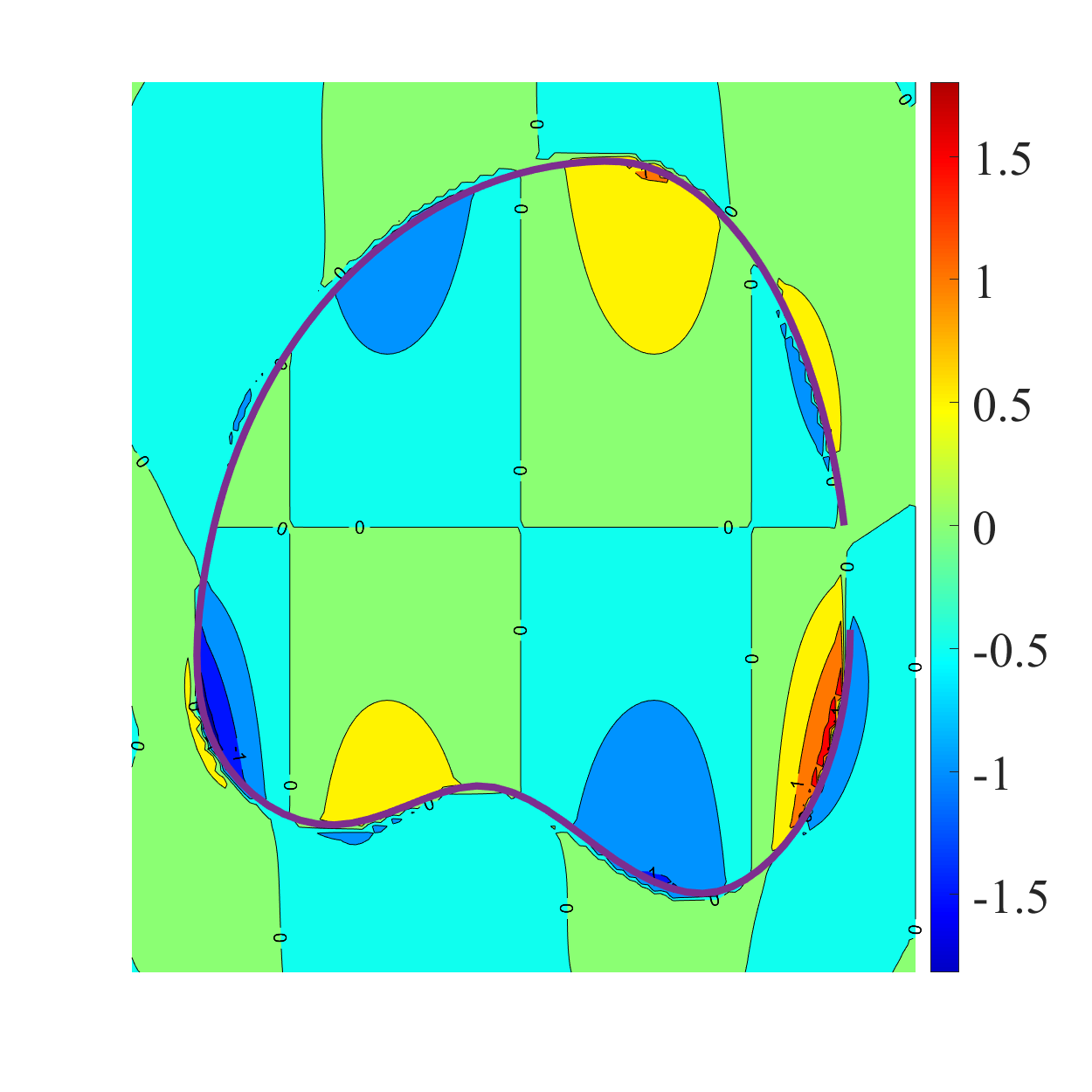}}
    \subfigure[The pressure field $p$]{\includegraphics[width=0.3\textwidth]{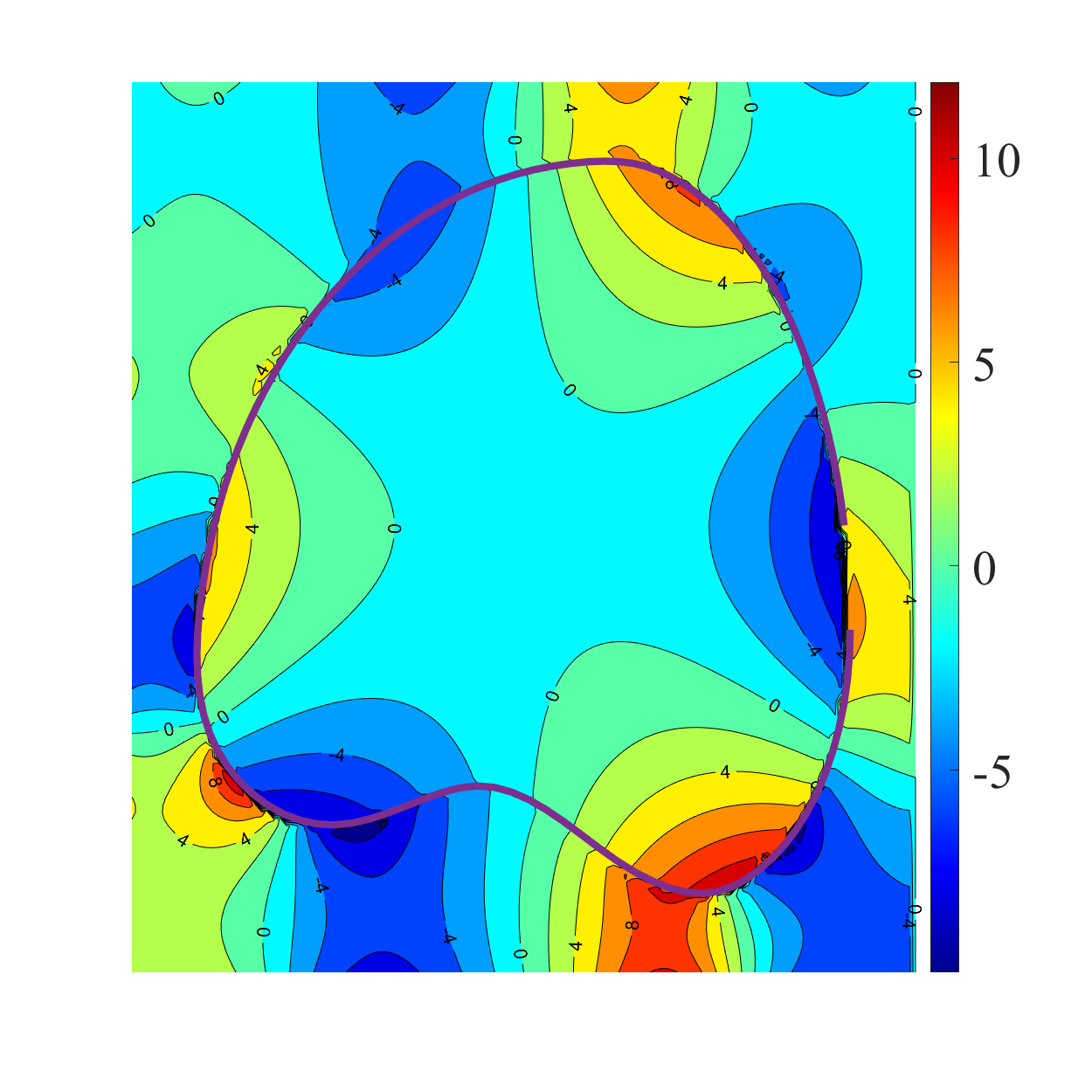}}
    \caption{The numerical solutions  for example 2 on the $512 \times 512$ grid.}
    \label{Stokes2D}
\end{figure}

\textbf{Example 2.}
This example solves the boundary value problem of the Stokes equation on the heart-shaped domain, with the Dirichlet boundary condition. A Cartesian grid-based MAC Scheme is applied to solve the Stokes equation. This approach places the pressure $p$ at the cell center, the $x-$component velocity $u^{(1)}$ and the $y-$component velocity $u^{(2)}$ at the midpoints of the vertical and horizontal edges of each cell, respectively. The method is detailed in\cite{dong2023second}. The boundary conditions are chosen so that the exact solution reads
\begin{equation}
\begin{aligned}
&u^{(1)}(x, y) = x(x^2 - 3y^2) + 1.5(1 - (x^2 + y^2))x\\
& u^{(2)}(x, y) = - 1.5(1 - (x^2 + y^2)^2)y\\
&  p(x, y) = 6(y^2 - x^2)\\
\end{aligned}
\end{equation}
 The bounding box $\mathcal{B}$ for the interface problem is set to be $\mathcal{B}=(-1.2,1.2) \times(-1.2,1.2)$. The execution time on the CPU and GPU are summarized in Tab.\,\ref{tab:Stokes}. Numerical results  are plotted in Fig.\,\ref{Stokes2D}. 
\begin{table}[H]
\centering
\begin{tabular}{c|c|ccccc}
\hline Boundary & \text { grid size }  & $64^2$ & $128^2$ & $256^2$ &  $512^2$ \\
\hline Dirichlet & \text {CPU time}& $6.95 \mathrm{~s}$ & $31.27 \mathrm{~s}$ & $135.26 \mathrm{~s}$ & $551.24\mathrm{~s}$\\
& \text{GPU time}&$1.36 \mathrm{~s}$ & $1.85 \mathrm{~s}$ & $3.25 \mathrm{~s}$ & $5.51\mathrm{~s}$   \\
& \text{Speedup}&$5.11$ & $16.90 $ & $41.61$  & $100.04$ \\
\hline
\end{tabular}
    \caption{BVP of the Stokes equation on the heart-shaped domain.}
    \label{tab:Stokes}
\end{table}

\textbf{Example 3.}
 This example solves the Gray-Scott model which consists of two singularly perturbed reaction-diffusion equations given by
$$
\begin{aligned}
&u_{t}=\epsilon_{1} \Delta u+\frac{1}{\epsilon_{0}}\left[\gamma(1-u)-u v^{2}\right] \\
&v_{t}=\epsilon_{2} \Delta v+\frac{1}{\epsilon_{0}}\left[u v^{2}-(\gamma+\kappa) v\right]
\end{aligned}
$$

Here, $u=u(x, y, t)$ and $v=v(x, y, t)$ are two unknown smooth functions, describing the concentration of some chemical substances in a bounded domain $\Omega$ for $t \geqslant 0 ; \gamma$ and $\kappa$ are respectively the feed and removal rate; $\epsilon_{0}, \epsilon_{1}$ and $\epsilon_{2}$ are small reactive or diffusive coefficients. In this example, the model is assumed to satisfy the homogeneous Neumann boundary condition that $\partial_{\mathrm{n}} u=\partial_{\mathrm{n}} v=0$ on $\partial \Omega$, initial condition and the involved parameters are specified as follows
$$
\begin{aligned}
&v(x, y, 0)=\left\{\begin{array}{cc}
\frac{1}{4} \sin ^{2}(4 \pi x) \sin ^{2}(4 \pi y), & -0.25 \leqslant x, y \leqslant 0.25 \\
0, & \text { others. }
\end{array}\right. \\
&u(x, y, 0)=1-2 v(x, y, 0) \\
&\gamma=0.024, \kappa=0.06, \epsilon_{0}=0.01, \epsilon_{1}=0.008, \epsilon_{2}=0.004
\end{aligned}
$$

\begin{figure}[ht]
    \centering{\includegraphics[width=1.0\textwidth]{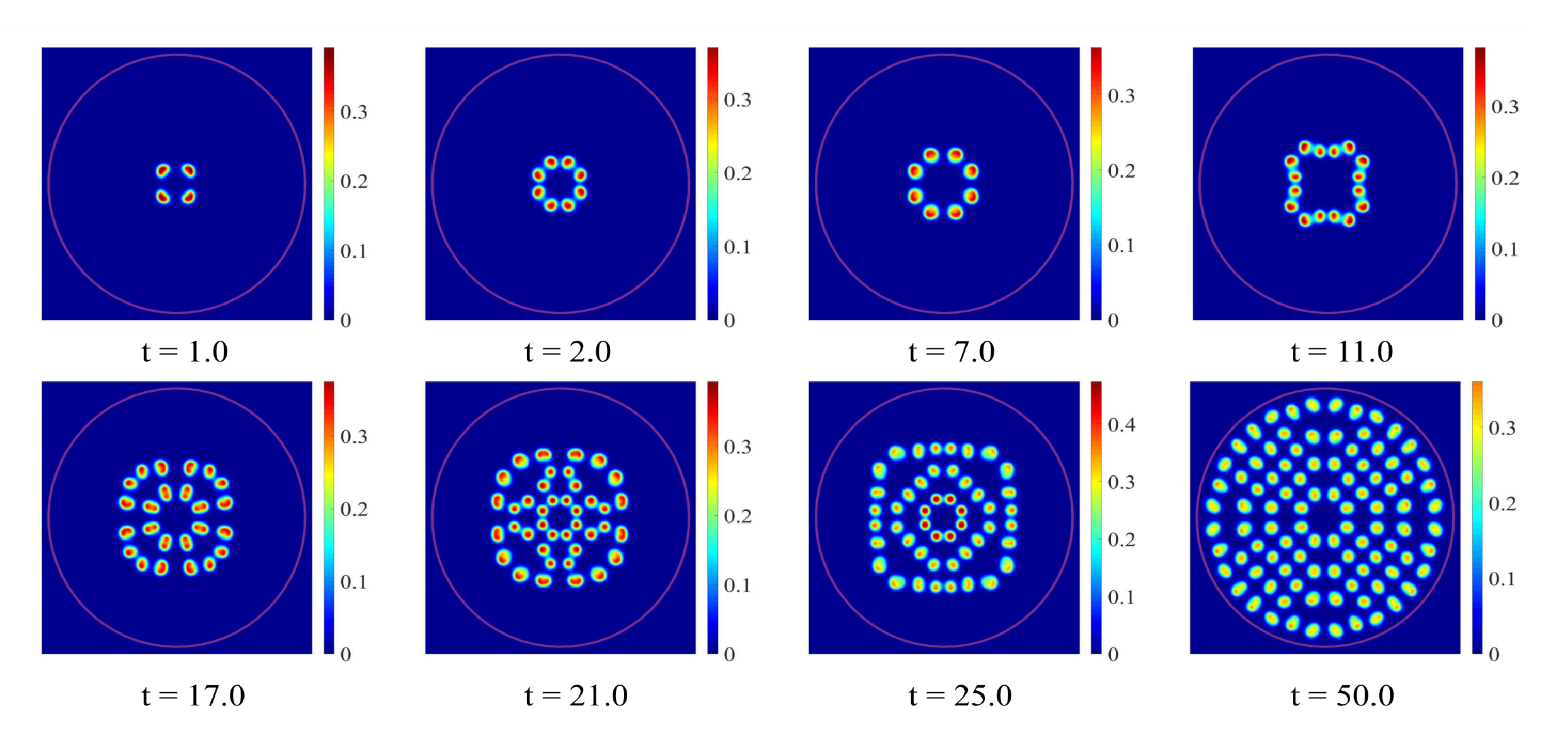}}
    \caption{the radius $r = 1.8$, $T = 1,2,7,11,17,21,25,50$ on the $128 \times 128$ grid.The bounding box $\mathcal{B}$  is set to be $\mathcal{B}=(-2.0,2.0) \times(-2.0,2.0)$ .}
    \label{PKFBIRD}
\end{figure}

The bounding box $\mathcal{B}$ for the interface problem is set to be $\mathcal{B}=(-2.0,2.0) \times(-2.0,2.0)$ and the tolerance is $10^{-8}$. Time direction is discretized by the second-order Strang splitting method\cite{MacNamara2016}. The numerical results when $T = 1,2,7,11,17,21,25,50$ are plotted in Fig.\,\ref{PKFBIRD}. Tab.\,\ref{tab:PKFBIRD} we present the execution time on the CPU and GPU of the parallel algorithm for different 
computing scales, In order to verify the computational efficiency and stability, the GPU acceleration ratio and numerical accuracy are also shown in the table. It is calculated selecting four different problem sizes: $128 \times 128$,  $256 \times 256$, $512 \times 512$, $1024 \times 1024$, the time step is increasing with the increase of the grid size. From the table we can see that the GPU acceleration ratio increases with increasing of the computation scale, It can be seen that a better performance is obtained when large problems are considered, which means our parallel method scales well.
\begin{table}[htpt!]

\centering
\begin{tabular}{c|ccccc}
\hline  boundary condition & \text { grid size }  & $128 \times 128$ & $256 \times 256$ & $512 \times 512$ &  $1024 \times 1024$\\
  &\text {Time steps}& $ 8$ & 16 & 32 &  64 \\
\hline Neumann & \text {CPU time}& $1.21 \mathrm{~s}$ & $ 8.98 \mathrm{~s}$ & $74.97 \mathrm{~s}$ & $633.78\mathrm{~s}$ \\

 & GPU time & $0.40 \mathrm{~s}$ & $1.04 \mathrm{~s}$ & $3.29 \mathrm{~s}$ & $12.07\mathrm{~s}$  \\
  & Speedup & $3.0250$ & $8.6346$ & $22.7872$ & $52.5086
$  \\
\hline
\end{tabular}
            \caption{Simulation time of CPU-based and GPU-based KFBI method, as well as the speedup  achieved by  GPU-based
 solver over the CPU-based solver. }
            \label{tab:PKFBIRD}
\end{table}

\iffalse
% \textbf{Example 4.Laplace equation in 3D torus domain}
% \begin{table}[H]
% \centering
% \begin{tabular}{c|c|ccccc}
% \hline Boundary & \text { grid size }  & $32^3$ & $64^3$ & $128^3$ &  $256^3$ & $512^3$\\
% \hline & \text {CPU time}& $0.64 \mathrm{~s}$ & $3.03 \mathrm{~s}$ & $21.87 \mathrm{~s}$ & $158.91\mathrm{~s}$ & $1654.21\mathrm{~s}$\\
% Dirichlet & GPU time&$0.56 \mathrm{~s}$ & $0.59 \mathrm{~s}$ & $0.68 \mathrm{~s}$ & $1.26\mathrm{~s}$  & $3.98 \mathrm{~s}$ \\
% &\left\|\mathbf{e}_{h}\right\|_{\infty} & 7.5 \mathrm{E}-4 & 1.0 \mathrm{E}-4 & 1.6 \mathrm{E}-5 &2.1  \mathrm{E}-6&3.2  \mathrm{E}-7 \\
% \hline
% \end{tabular}
%     \caption{\footnotesize{DBVP of the Laplace equation in torus domain with  bounding box $\mathcal{B}=(0,1.2) \times(0,1.2)\times(0,1.2)$.}}
%     \label{tab:my_label}
% \end{table}

% \begin{figure}[H]
%     \centering{\includegraphics[width=0.45\textwidth]{figure/3D.png}}
%     \caption{\footnotesize{The numerical solutions with the $256 \times 256 \times 256$ grid.}}
% \end{figure}
\fi
\textbf{Example 4.}
 This example solves the Dirichlet BVP of the Stokes equation on an sphere $\Omega$ which is given by
\begin{equation}
\Omega=\left\{(x, y, z) \in \mathbb{R}^3: \frac{x^2}{a^2}+\frac{y^2}{b^2}+\frac{z^2}{c^2}<1\right\}
\label{sphere}
\end{equation}
with $a=1.0, b=0.8, c=0.6$.  The bounding box $\mathcal{B}$ for the interface problem is $\mathcal{B}=[-1.2,1.2] \times[-1.2,1.2] \times[-1.2,1.2]$.   The Dirichlet BC is chosen so that the exact solution reads
% \begin{equation}
% \begin{aligned}
% & u^{(1)}(x, y, z)= \begin{cases}\exp (\cos y)+\exp (\sin z), & x^2+y^2+z^2>1, \\
% -4\left(1-x^2-y^2\right) x y-4 x^2 z^2+\left(x^2+3 z^2-2\right)\left(z^2-x^2\right), & x^2+y^2+z^2 \leq 1,\end{cases} \\
% & u^{(2)}(x, y, z)= \begin{cases}\exp (\sin x), & x^2+y^2+z^2>1, \\
% -4 x^2 y^2+\left(3 x^2+y^2-2\right)\left(x^2-y^2\right), & x^2+y^2+z^2 \leq 1,\end{cases} \\
% & u^{(3)}(x, y, z)= \begin{cases}\exp (\cos (x)), & x^2+y^2+z^2>1, \\
% -4\left(1-x^2-z^2\right) x z, & x^2+y^2+z^2 \leq 1\end{cases} \\
% & p(x, y, z)= \begin{cases}\exp (\cos x+\sin y)+\exp (\cos z+\sin x), & x^2+y^2+z^2>1, \\
% (x-1)^3+(y-1)^3+(z-1)^2, & x^2+y^2+z^2 \leq 1 .\end{cases} \\
% &
% \end{aligned}
% \end{equation}
\begin{equation}
\begin{aligned}
& u^{(1)}(x, y, z)= \exp (\cos y)+\exp (\sin z)-4\left(1-x^2-y^2\right) x y-4 x^2 z^2+\left(x^2+3 z^2-2\right)\left(z^2-x^2\right) \\
& u^{(2)}(x, y, z)= \exp (\sin x)-4 x^2 y^2+\left(3 x^2+y^2-2\right)\left(x^2-y^2\right)\\
& u^{(3)}(x, y, z)= \exp (\cos (x))-4\left(1-x^2-z^2\right) x z \\
& p(x, y, z)= \exp (\cos x+\sin y)+\exp (\cos z+\sin x)+8  (3  x^2 - y^2)  y + 8  x  (3  z^2 - x^2) \\
&
\end{aligned}
\end{equation}

The error orders and execution times on the CPU and GPU are encapsulated in Table \ref{tab:Stokes3D}, derived from four distinct problem sizes: $32^3$,  $64^3$, $128^3$, $256^3$, and $512^3$. The table reveals an increasing GPU acceleration ratio with the escalation of computational scale. It is observed that enhanced performance is achieved for larger problems, indicating the scalability of the parallel method.
\begin{table}[H]
\centering
\begin{tabular}{c|c|ccccc}
\hline Boundary & \text { grid size }  & $32^3$ & $64^3$ & $128^3$ &  $256^3$ \\
\hline & \text {CPU time} & $44.85 \mathrm{~s}$ & $272.39 \mathrm{~s}$ & $1948.24\mathrm{~s}$ &$13521.36\mathrm{~s}$\\
Dirichlet & \text{GPU time} & $1.34 \mathrm{~s}$ & $2.82 \mathrm{~s}$ & $4.98\mathrm{~s}$  & $38.62 \mathrm{~s}$ \\
 & \text{Speedup} & $33.47$ & $96.59$ & $391.21$  & $350.11$ \\
%&\left\|\mathbf{e}_{h}\right\|_{\infty} &  5.69 \mathrm{E}-2 & 5.89 \mathrm{E}-3 &8.78  \mathrm{E}-4&1.04 \mathrm{E}-4 \\
\hline
\end{tabular}
    \caption{BVP of the Stokes equation  on the bounding box $\mathcal{B}=(-1.2,1.2) \times(-1.2,1.2)\times(-1.2,1.2)$ with GMRES iterition method.}
    \label{tab:Stokes3D}
\end{table}

% \begin{figure}[htb]
%     \centering
%     \subfigure[The velocity field $u^{(1)}$]{\includegraphics[width=0.4\textwidth]{figure/u1.eps}}
%     \subfigure[The velocity field $u^{(2)}$]{\includegraphics[width=0.4\textwidth]{figure/u2.eps}}
%     \subfigure[The velocity field $u^{(3)}$]{\includegraphics[width=0.4\textwidth]{figure/u3.eps}}
%     \subfigure[The pressure field $p$]{\includegraphics[width=0.4\textwidth]{figure/p.eps}}
%     \caption{The numerical solutions  for example 5 on the $128 \times 128 \times 128$ grid.}
%     \label{Stokes3D}
% \end{figure}
\begin{figure}[htpt!]
    \centering{\includegraphics[width=1.05\textwidth]{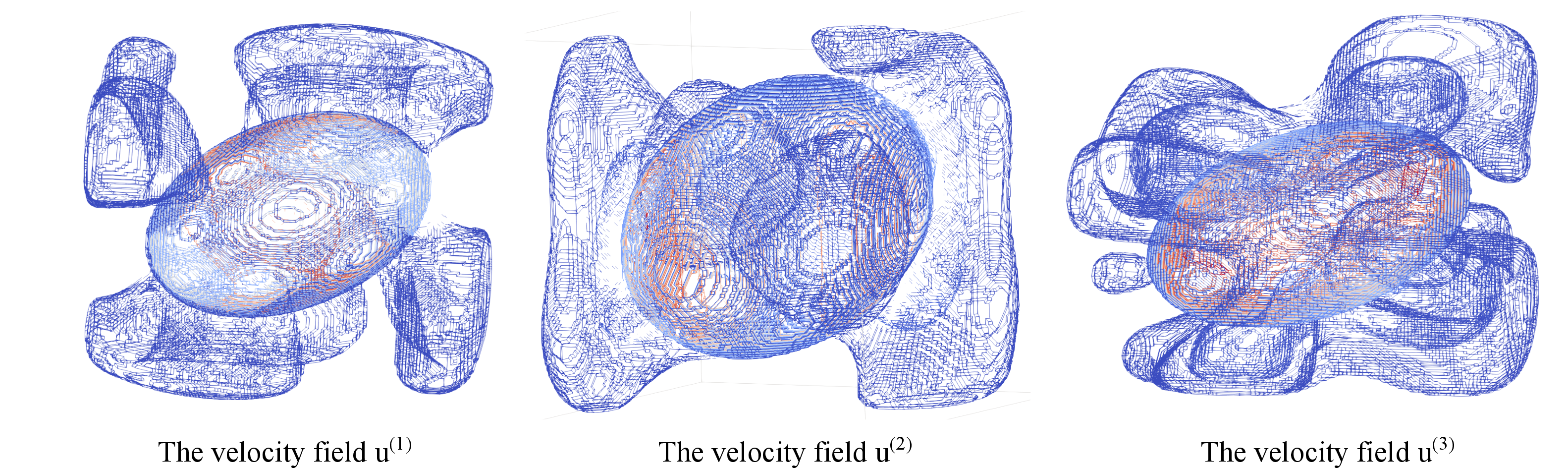}}
    %\caption{The numerical solutions for feature edge in example 5 on the $128 \times 128 \times 128$ grid.}
    \label{One_Poisson}
\end{figure}

\subsection{Multiple-GPU results}
To augment computational precision, refinement
is applied to the computational grid. Example 6 uses the same numerical test cases as examples 1, 3, and 4. Table \ref{tab:multi_gpu} displays the computation times for solving the 2D Laplace, reaction-diffusion, and 3D Stokes equation. These computations were executed using 1, 2, 4, and 8 GPUs.

\begin{figure}[htb]
    \centering
    \subfigure[2D Laplace]{\includegraphics[width=0.31\textwidth]{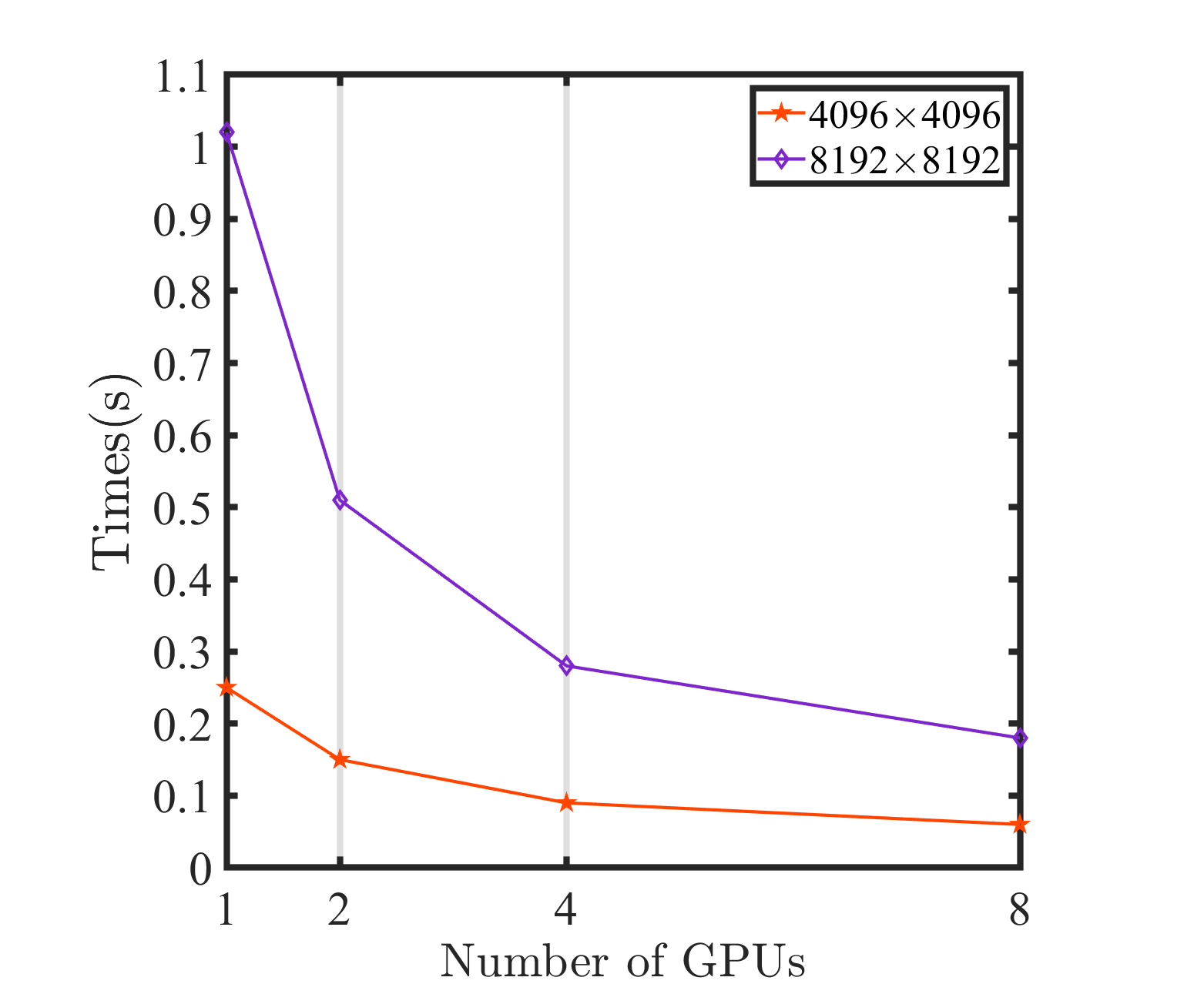}}
    \subfigure[2D reaction-diffusion]{\includegraphics[width=0.31\textwidth]{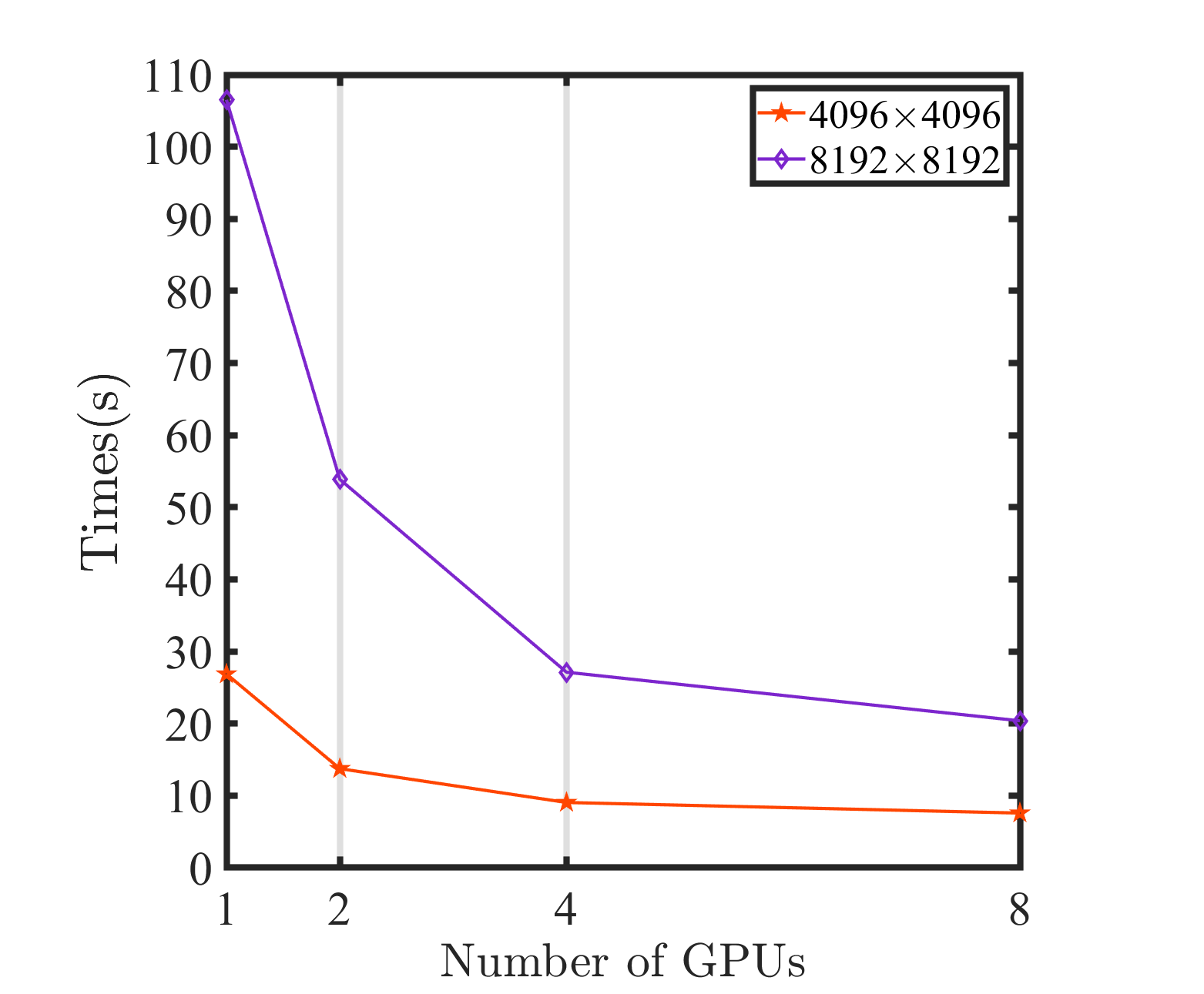}}
    \subfigure[3D Stokes]{\includegraphics[width=0.31\textwidth]{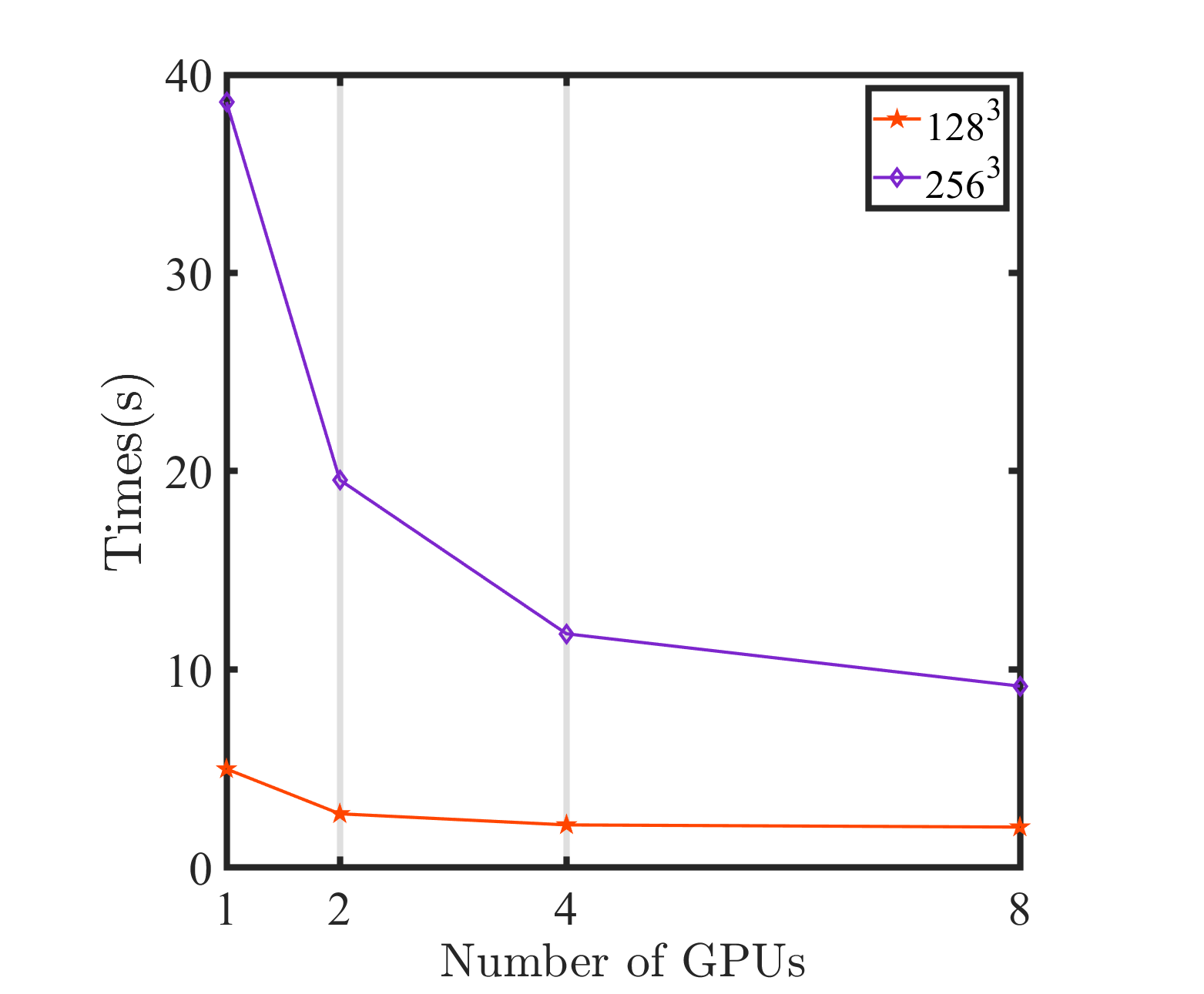}}
    \caption{Comparison of parallel efficiency with different numbers of GPUs}
    \label{MULTIGPU}
\end{figure}

In the Fig.\,\ref{MULTIGPU}, we can conclude that multi-GPUs parallel computing achieves linear speedup, despite a slight decrease in single-GPU performance when more GPUs are employed, as shown in Tab.\,\ref{tab:multi_gpu}. The linear growth of parallel efficiency hindrance can be attributed to inter-GPU communication, involving tasks such as the distribution of boundary data(point 1 of procedure 2), exchange of ghost cell data(points 3 and 4 of procedure 2), and collection of boundary data(the point 6 of procedure 2) in \ref{Sec:Algorthm summary}.

\begin{table}[H]
\centering
\footnotesize
\scalebox{0.8}{
\begin{tabular}{c|ccc|ccc|ccc}
\hline \text { equation } & & 2D Laplace && &2D reaction-diffusion &&& 3D Stokes&\\
\hline  \text {grid size }  & $4096^2$ & $8192^2$ & $16384^2$ & $4096^2$ & $8192^2$ & $16384^2$&$128^3$ &  $256^3$ &  $512^3$ \\
\hline $1$GPU & $0.25 \mathrm{~s}$ & $1.02 \mathrm{~s}$ &  &  $26.87 \mathrm{~s}$ & $106.51 \mathrm{~s}$   &  & $4.98\mathrm{~s}$ &$38.62\mathrm{~s}$ \\
   $2$GPUs  & $0.15 \mathrm{~s}$ & $0.51 \mathrm{~s}$&& $13.74 \mathrm{~s}$ & $53.87 \mathrm{~s}$ & & $2.72 \mathrm{~s}$ & $19.56 \mathrm{~s}$\\
    $4$GPUs   & $0.09 \mathrm{~s}$ & $0.28 \mathrm{~s}$&$1.01 \mathrm{~s}$& $9.05 \mathrm{~s}$ & $27.12 \mathrm{~s}$ &$123.78 \mathrm{~s}$  &  $2.16 \mathrm{~s}$ &  $11.8 \mathrm{~s}$  \\
      $8$GPUs   & $0.06 \mathrm{~s}$ & $0.18 \mathrm{~s}$&$0.62 \mathrm{~s}$& $7.58 \mathrm{~s}$ & $20.38 \mathrm{~s}$ & $74.93\mathrm{~s}$ &  $2.05 \mathrm{~s}$ &  $9.16 \mathrm{~s}$ &$66.21 \mathrm{~s}$ \\
      
\hline
\end{tabular}
}
    \caption{Multi-GPUs execution time vs. single GPU.}
    \label{tab:multi_gpu}
\end{table}

\section{Discussion}

This paper presents a second-order, single, and multiple-GPU accelerated efficient KFBI method for  elliptic boundary value problems. The equations are first transformed into a BIE, and then the unknown density in the equation is solved by GMRES iteration. Boundary and volume integral can be evaluated by equivalent interface problems to obtain the approximate solution. The procedure for solving the interface problem consists of four steps: discretization, correction, fast solving, and interpolation.
 
In the single GPU algorithm, since the KFBI method itself mainly 
focuses on the control points on the boundary and the irregular nodes near the boundary. We only need to assign  threads to them and design a fast algorithm on the GPU to solve the interface problem efficiently. In the multiple-GPU algorithm, the system of linear equations in the interface problem must be solved by the ADM method, which involves the interaction of information between the host and the device handled by MPI. 

The accuracy and efficiency of the algorithm are verified from numerical examples. The method is especially suited for GPU acceleration in large-scale computations, and the multiple-GPU distributed solver scales well. Numerical examples show that single-GPU solver speeds 50-200 times than traditional CPU while the eight GPUs distributed solver yields up to $60\%$ parallel efficiency.

The single/multiple-GPU accelerated KFBI method can be extended for other PDEs, such as the Maxwell and elasticity equations. Furthermore, combined with the deep learning method, the KFBI method may exhibit potential applicability in solving equations within non-smooth domains on the GPU platform.
\section*{Acknowledgments}
This work is financially supported by the Strategic Priority Research Program of Chinese Academy of Sciences(Grant No. XDA25010405). It is also partially  supported by the National Key R\&D Program of China, Project Number 2020YFA0712000, the National Natural Science Foundation of China (Grant No. DMS-11771290) and the Science Challenge Project of China (Grant No. TZ2016002). Additionally, it is supported by the Fundamental Research Funds for the Central Universities. 

%%%% Bibliography  %%%%%%%%%%
% \begin{thebibliography}{99}
% \bibitem{Berger}M. J. Berger and P. Collela, Local adaptive mesh refinement
% for shock hydrodynamics,
% J. Comput. Phys., 82 (1989), 62-84.
% \bibitem{deBoor}C. de Boor,  Good Approximation By Splines With Variable Knots II, in Springer Lecture
%  Notes Series 363, Springer-Verlag, Berlin, 1973.
% \bibitem{TanTZ} Z. J. Tan, T. Tang and Z. R. Zhang, A simple moving mesh method for one- and
% two-dimensional phase-field equations, J. Comput. Appl. Math., to appear.
% \bibitem{Toro}E. F. Toro, Riemann Solvers and Numerical Methods for Fluid Dynamics,
% Springer-Verlag Berlin Heidelbert, 1999.
% \end{thebibliography}

\bibliographystyle{unsrt} 
\bibliography{ref} 
\end{document}